\documentclass[12pt,twoside,leqno]{amsart}
\usepackage[cmtip, arrow]{xy}
\usepackage{amsmath,amscd,amssymb,amsfonts,pb-diagram,pb-xy}
\usepackage{graphicx}

\newcommand{\Z}{\mathbb{Z}}
\newcommand{\Q}{\mathbb{Q}}

\newcommand{\C}{\mathbb{C}}
\newcommand{\F}{\mathbb{F}}
\newcommand{\V}{\mathbb{V}}

\renewcommand{\H}{\mathbb{H}}
\renewcommand{\P}{\mathbb{P}}

\newcommand{\mO}{\mathcal{O}}
\newcommand{\mP}{\mathcal{P}}
\newcommand{\mQ}{\mathcal{Q}}
\newcommand{\mY}{\mathcal{Y}}
\newcommand{\mZ}{\mathcal{Z}}

\newcommand{\lra}{\longrightarrow}

\theoremstyle{plain}

\numberwithin{thm}{section}
\numberwithin{equation}{section}

\begin{document}
\title
{Picard-Fuchs operators\\
for octic arrangements I\\
(\em The case of orphans)\\
}

\author{Slawomir Cynk and Duco van Straten}
\date{}
\thanks{The first author was partly supported by the Schwerpunkt Polen (Mainz) and NCN grant no. N N201 608040. 
The second author was partly supported by DFG Sonderforschungsbereich/Transregio 45.}
\begin{abstract}
We report on $25$ families of projective Calabi-Yau threefolds that do
not have a point of maximal unipotent monodromy in their moduli space. 
The construction is based on an analysis of certain pencils of octic 
arrangements that were found by {\sc C. Meyer} \cite{Mey}. 
There are seven cases where the Picard-Fuchs operator is of order two and $18$ 
cases where it is of order  four. The birational nature of the Picard-Fuchs equation can be used  effectively to distinguish between families whose members have the same Hodge numbers. 
\end{abstract}

\maketitle
\section{Introduction}
The phenomenon of mirror symmetry among Calabi-Yau threefolds
has attracted a lot of attention and has led to major developments in mathematics
and physics, see e.g. \cite{Horics}, \cite{GHJ}. Especially the marvelous discovery 
by {\sc P. Candelas}, {\sc X. de la Ossa} and coworkers \cite{COGP} of the relation 
between enumeration of rational curves on a Calabi-Yau threefold and period integrals 
on another mirror manifold has been an inspiration to many researchers. Both the 
determination  of instanton numbers in \cite{COGP} and the construction of mirror 
pairs by {\sc V. Batyrev} \cite{Bat1} as generalised by {\sc M. Gross} and 
{\sc B. Siebert} \cite{GS}, depend on {\em families} of Calabi-Yau manifolds that 
degenerate at the boundary of their moduli space at a point of {\em maximal unipotent monodromy} (MUM), \cite{Mor}. In many natural families of Calabi-Yau threefolds, like hypersurfaces in toric varieties defined by reflexive polytopes, there do exist 
such MUM-points in their moduli space.\\

However, it has been known for some time that there are quite simple
examples of one-parameter families 
\[ f: \mathcal{Y} \lra S\]
of Calabi-Yau threefolds  for which there are no such MUM points.
First examples of this kind were described by {\sc J.-C. Rohde} \cite{Roh1} 
and later by {\sc A. Garbagnati} and {\sc B. van Geemen} \cite{GG}. 
However, these examples are somewhat atypical in the sense that the 
cohomological local system 
$\H=R^3f_*(\C)$ decomposes as a tensor product 
\[ \H= \F \otimes \V ,\]
where $\F$ is a constant VHS with Hodge numbers $(1,0,1)$ and $\V$ a variable 
$(1,1)$-VHS. As a consequence, although the cohomology space is four
dimensional, the Picard-Fuchs operator is of order two, and {\sc
  Rohde} \cite{Roh2} asked  
the question if there exist examples where the Picard-Fuchs operator had order 
four.

In \cite{Str} we described an example of such a family, but the members of
that family had the defect of not being projective. {\sc W. Zudilin} \cite{Zud} 
suggested to name {\em orphan} to describe such families, as they do not have a MUM. 
We decided to follow his suggestion, and in this paper we will describe a series of new orphans which have the virtue of being projective as well. The first such example  was announced in \cite{CvS1} and has the full symplectic 
group $Sp_4(\C)$ as differential Galois group.\\ 

\centerline{\bf Main Result}
\vskip 10pt
{\em There exist at least $10$ birationally unrelated families of projective orphans with Picard-Fuchs operator of order four.}\\

For a precise definition of the notion of {\em related families} we refer to the
discussion in section $7$.\\

We also want to point out two interesting phenomena that we discovered during
the analysis of the examples.
\begin{itemize}
\item We found one new (conjectural) instance of a Calabi-Yau threefold
with a Hilbert modular form of weight $(4,2)$ and level $6\sqrt{2}$, much
like the famous example of {\sc Consani} and {\sc Scholten}, \cite{CS}.
(Hilbert modularity for that example was shown in \cite{DPS}.)

\item We found an example of cohomology change in a family of Calabi-Yau
threefolds at a point with monodromy of finite order. As a consequence, the
central fibre of any semi-stable model is reducible, contrary to what happens
in families of K3-surfaces, according to {\sc Kulikov}s theorem.  
\end{itemize}

\section{Double octics}
By a {\em double octic} we understand a double cover $Y$ of $\P^3$ 
ramified over a surface of degree $8$. It can be given by an equation 
of the form
\[ u^2=f_8(x,y,z,w)\]
and thus can be seen as a hypersurface in weighted projective space
$\P(1,1,1,1,4)$. For a general choice of the degree eight polynomial $f_8$ 
the variety $Y$ is a smooth Calabi-Yau space with Hodge numbers 
$h^{11}=1, h^{12}=149$. When the octic $f_8$ is a product of eight linear 
factors, we speak of an {\em octic arrangement}. These form a nine dimensional
sub-family and for these, the double cover $Y$ has singularities at the 
intersections of the planes. In the generic such situation $Y$ is singular 
along  $8.7/2=28$ lines, and by blowing up these lines (in any order) we 
obtain a smooth Calabi-Yau manifold $\widetilde{Y}$ with $h^{11}=29, h^{12}=9$. 

By taking the eight planes in special positions, the double cover $Y$ acquires 
further singularities. As explained in \cite{Mey}, if the arrangement does not 
have double planes, fourfold lines or sixfold points, there exist a diagram
\[ 
\begin{diagram}
\node{\widehat{Y}} \arrow{e,t}{\hat{\pi}}  \arrow{s,l}{2:1} \node{Y}\arrow{s,r}{2:1}\\
\node{\widehat{\P}^3} \arrow{e,t}{\pi} \node{\P^3}
\end{diagram}
\]
where $\pi$ is a sequence of blow ups, $\hat{\pi}$ a crepant resolution, and the
vertical maps are two-fold ramified covers.\\

In this way a myriad of different Calabi-Yau threefolds $\widehat{Y}$ can be 
constructed. One of the nice things is that one can read off the Hodge number $h^{12}$
as the dimension of the space of deformations of the arrangements that do not change
the combinatorial type, \cite{CvS}. In his doctoral thesis \cite{Mey}, 
{\sc C. Meyer} found $450$ combinatorially different octic arrangements, determined their Hodge numbers and started the study of  their arithmetical properties.\\
 
Among these $450$ arrangements there were $11$ arrangements with 
$h^{12}(\widehat{Y})=0$, so lead to rigid Calabi-Yau threefolds and $63$ one-parameter families of arrangements with $h^{12}(\widehat{Y}_t)$ leading to one-parameter families (all defined over $\Q$) of Calabi-Yau threefolds, parametrised by
$\P^1$: for general $t \in \P^1$, the crepant resolutions
\[ 
\begin{diagram}
\node{\widehat{Y}_t} \arrow{e,t}{\hat{\pi}}\node{Y_t}
\end{diagram}
\]

of the double octics $Y_t$ can be put together into a family $\overline{\mathcal{Y}}$ 
over $S:=\P^1\setminus \Sigma$, where $\Sigma \subset \P^1$ is a finite set of special values, where the combinatorial type
of the arrangement changes. At these points, the configuration becomes rigid, or ceases to be of Calabi-Yau type: the arrangement contains a double plane, a fourfold line or a sixfold point.
When we use the sequence of blow-ups to resolve the generic fibre and apply this to all members of the family, we arrive at a diagram of the form

\[ 
\begin{diagram}
\node{\widehat{\mathcal{Y}}} \arrow{e,J} \arrow{s,l}{f}\node{\widehat{\overline{\mathcal{Y}}}}\arrow{s,r}{\overline{f}}\\
\node{S} \arrow{e,J} \node{\P^1}
\end{diagram}
\]

where $f$ is a smooth map, with fibre $\widehat{Y}_t$, a smooth Calabi-Yau threefold with $h^{12}=1$.
In general, the fourfold $\hat{\overline{\mathcal{Y}}}$ will have singularities sitting over 
the special points $s \in \Sigma $.\\

Recently, in \cite{CKC} the analysis of {\sc Meyer} was found to be complete and only three further examples 
with $h^{12}=0$ exist over number fields and there is one further example of a family with $h^{12}=1$, defined over $\Q(\sqrt{-3})$. Furthermore, in that paper various symmetries and birational
maps between different arrangements were found.\\ 

We will take a closer look at these $63$ {\sc Meyer}-families. In order to facilitate comparison with literature,
we will keep the numbering from \cite{Mey}. In four cases (arrangements ${\bf 33}$, ${\bf 155}$,  
${\bf 275}$, ${\bf 276}$) small adjustments in the parametrisation of the family were made.\\

Although arithmetical information on varieties in these families is readily available via the counting 
of points in finite fields, we found it extremely hard to understand details of the resolution and topology 
from the combinatorics of the arrangements. For example, the Jordan type of the local monodromy around the 
special points turned out to be very delicate. In particular, we failed to find a clear combinatorial way to
recognise the appearance of a MUM-point.\\

As an example, we look at configuration ${\bf 69}$ of {\sc Meyer}. It consists of six planes making up a cube, with two additional planes that pass through a face-diagonal and opposite vertices of the cube, as in the
following picture.

\begin{center}
\includegraphics[height=5cm]{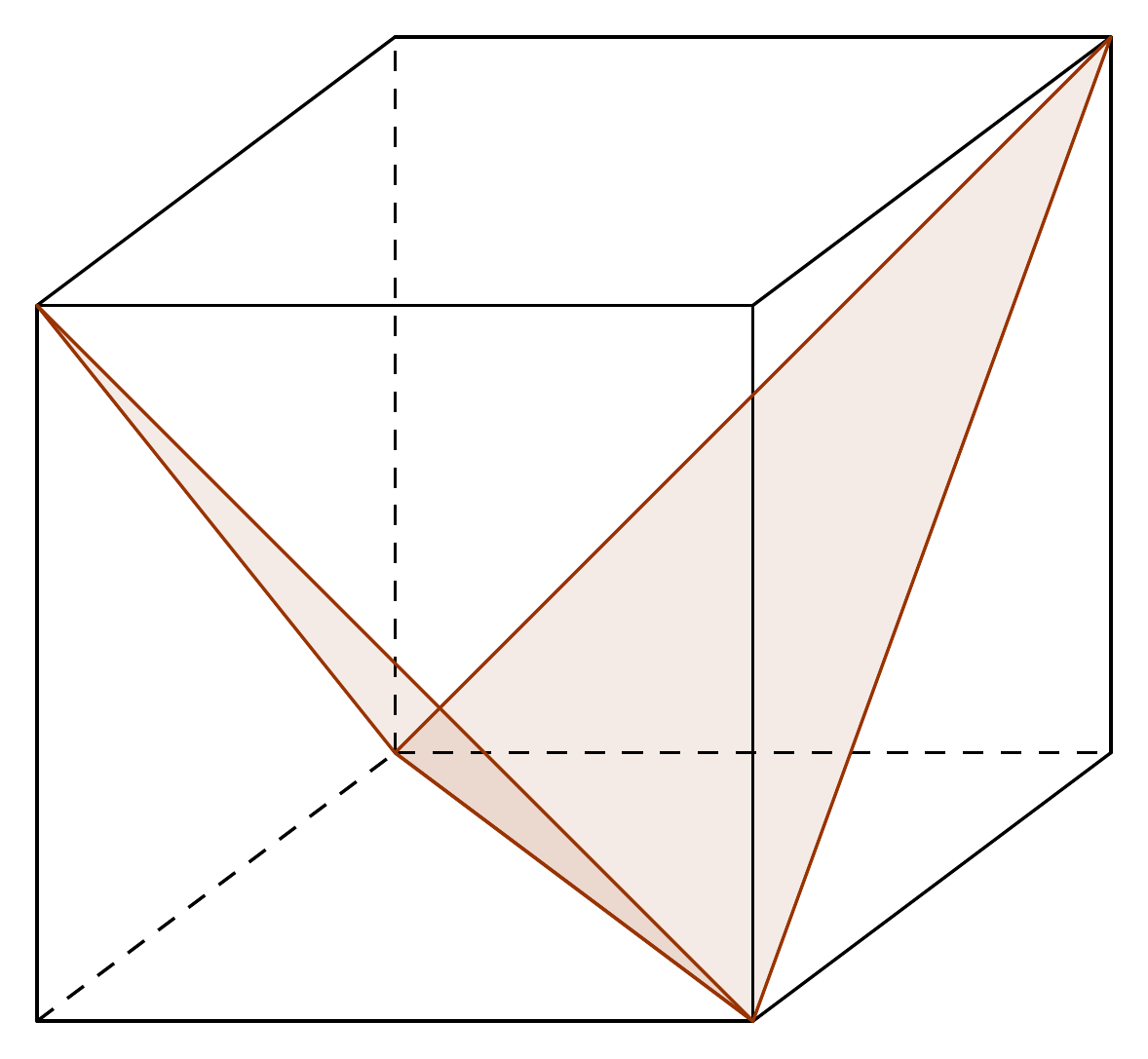}

{\bf \em Arrangement 69.}
\end{center}
The configuration is rigid and its resolution is a rigid Calabi-Yau with $h^{11}=50, h^{12}=0$. 
 
By sliding the intersection point at the corner of one of the two planes, we arrive at configuration 
${\bf 70}$, with $h^{11}=49, h^{12}=1$. Clearly, this pencil contains another rigid configuration, 
namely arrangement ${\bf 3}$.

\begin{center}
\includegraphics[height=5cm]{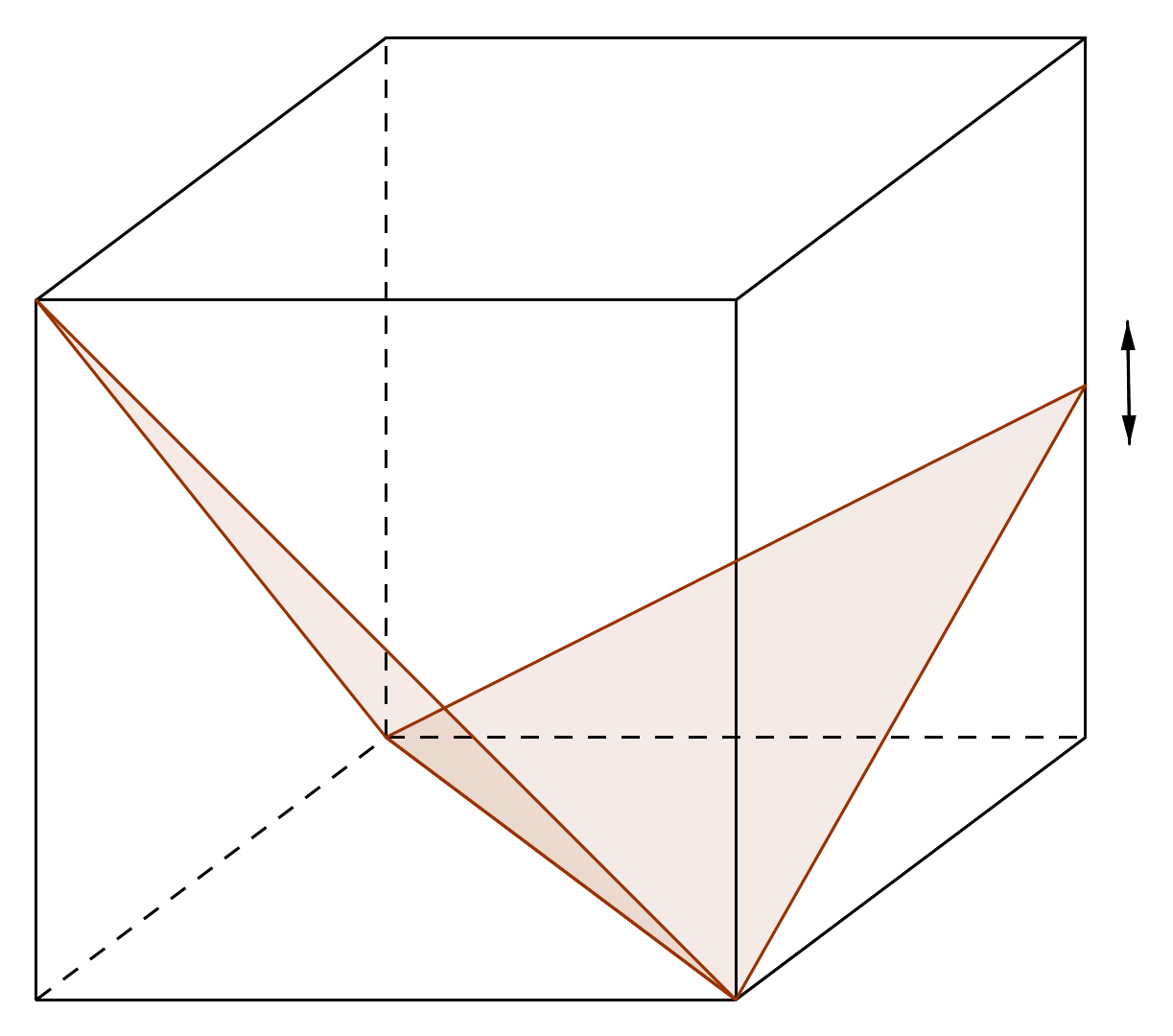}\\
\includegraphics[height=5cm]{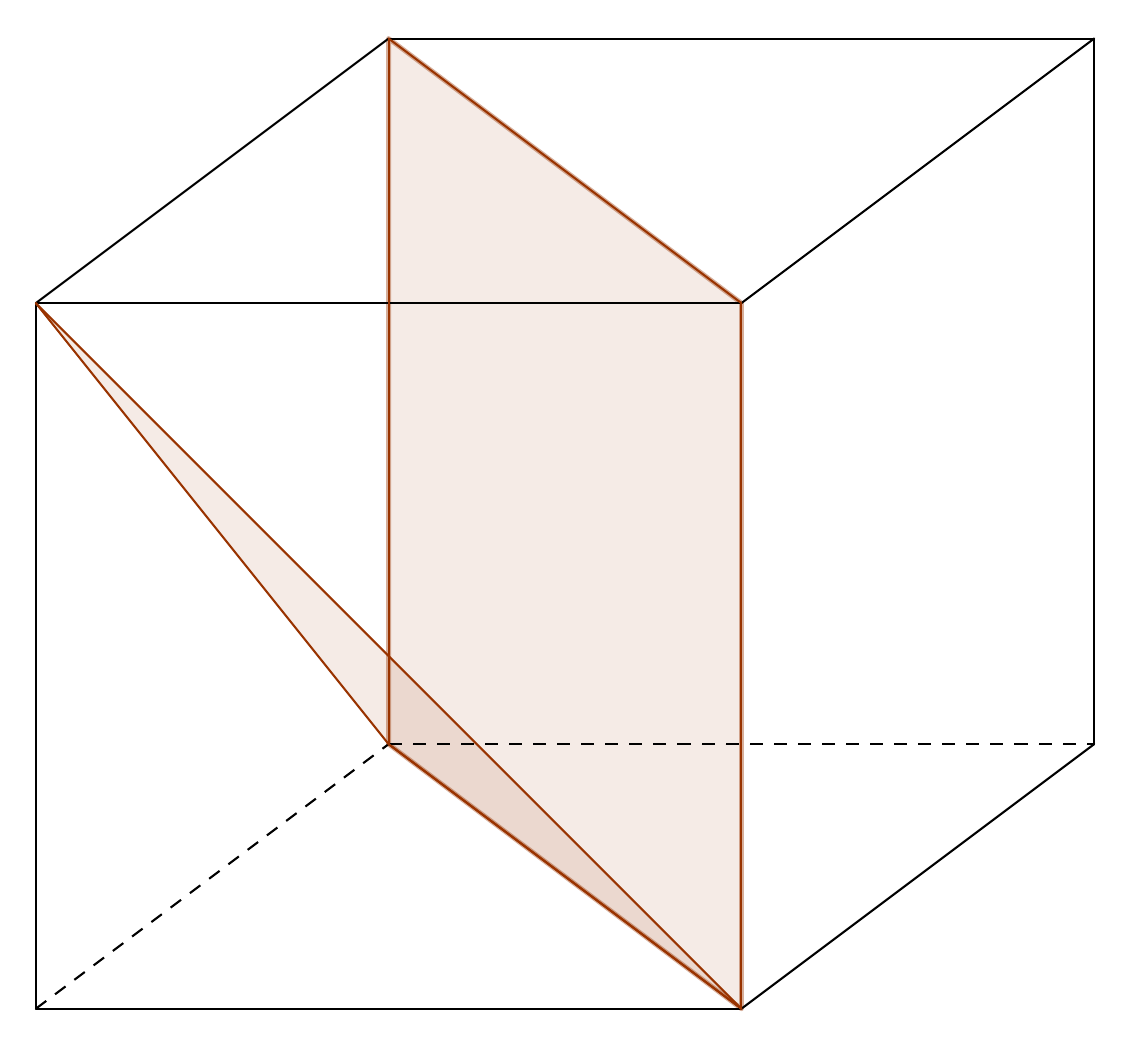}

{\bf \em Family of arrangements 70 and arrangement 3 }
\end{center}

But there are also two degenerations involving two double planes and it is not
so clear what the corresponding fibres of the semi-stable reduction look like, nor were we be able 
to determine topologically the monodromy around these points. It turns
out this family is one of the simplest orphans we know of and in this
paper we will be dealing with the $25$ cases without such a MUM-point
and which thus do not make it onto the list of Calabi-Yau operators
\cite{AZ}, \cite{AESZ}, \cite{CYDB}. In the sequel \cite{CvS2} to this
botanical paper, we will report on Picard-Fuchs operators for the remaining 
$38$ MUM-cases.

\section{$(1,1,1,1)$-variations}

\subsection{Generalities}

Let us consider more generally a fibre square diagram  
\[
\begin{diagram}
\node{\mathcal{Y}} \arrow{e,J} \arrow{s,l}{f} \node{\overline{\mathcal{Y}}} \arrow{s,r}{\overline{f}}\\
\node{S} \arrow{e,J} \node{\P^1}
\end{diagram}
\]

where $f: \mathcal{Y} \lra S:= \P^1 \setminus \Sigma$ is a smooth proper map of Calabi-Yau threefolds.
The datum of the local system
\[ \H_{\C}:=R^3f_* \C_{\mathcal{Y}}\]
is equivalent, after choice of a base point $s \in S$, to that of its 
monodromy representation
\[ \rho: \pi_1(\P^1\setminus \Sigma, s) \lra Aut(H^3(Y_s,\C))\]
There is an underlying lattice bundle $\H_{\Z}$, coming from the integral cohomology and a
non-degenerate skew-symmetric intersection pairing, causing this representation to land in
\[Aut(H^3(Y_s,\Z)/torsion)=Sp_m(\Z),\;\;\; m= \dim H^3(Y_s).\]

Furthermore, $\H$ carries the structure of a {\em polarised variation of Hodge structures (VHS)}, 
meaning basically that the fibres $\H_t=H^3(Y_t,\C)$ of the local system carry a pure (polarised)
Hodge structure with Hodge numbers $(1,h^{12},h^{12},1)$.
It is a fundamental fact, proven by {\sc Schmid} \cite{Schm1}, that one may complement this VHS
defined on $\P^1 \setminus \Sigma$ by adding for each $s \in \Sigma$ a so called {\em Mixed Hodge Structure} 
(MHS) $(H_s,W_{\bullet}, F^{\bullet})$. Here $H_s$ is a $\Q$-vector space that can be identified with 
the sections of the  $\Q$-local system $\H_{\Q}$ over an arbitrary small slit disc centered at $s$. 
The local monodromy transformation $T:=T_s:H_s \to H_s$ at $s$ can be written as 
\[T=U S\]
where $U$ is unipotent and $S$ is semi-simple. The {\em monodromy logarithm}
\[ N= -\log U=(1-U)+\frac{1}{2}(1-U)^2+\frac{1}{3}(1-U)^3+\ldots\]
is nilpotent and determines a weight filtration $W_{\bullet}$ on $H_s$
which is characterised by the property that $N:W_k \to W_{k-2}$ and
\[ N^k: Gr^W_{d+k} H_s \stackrel{\simeq}{\lra} Gr^W_{d-k} H_s.\]
The Hodge filtration $F^{\bullet}$ in $H_s \otimes \C$ arises as limit from the 
Hodge filtration on the spaces $\H_t$, when $t \mapsto s$, and it is a fundamental 
fact that for each $s\in \Sigma$ the $F^{\bullet}$ defines a pure Hodge structure of 
weight $k$ on the graded pieces $Gr^W_k H_s$.\\

In the geometrical case {\sc Steenbrink} \cite{Ste} has constructed this mixed Hodge 
structure on $H_s$ using a semi-stable model

\[ \begin{diagram}
\node{D} \arrow{e,J} \arrow{s} \node{\mZ} \arrow{e} \arrow{s} \node{\overline{\mY}}\arrow{s,r}{\overline{f}}\\
\node{\{s\}}\arrow{e,J} \node{\Delta} \arrow{e} \node{\P^1}
\end{diagram}
\]

Here $\Delta$ is a small disc, $\Delta \to \P^1$ is a finite covering map, 
ramified over one of the $s \in \Sigma$, $\mZ$ is smooth and the fibre $D$ over $s$
is a (reduced) normal crossing divisor with components $D_i$ inside $\mZ$.
The complex of relative logarithmic differential forms
\[\Omega_{\mZ/\Delta}^{\bullet}(\log D)\]
can be used to describe the cohomology of the fibres and its extension to
$\Delta$. The complex comes with two
filtrations $F^{\bullet}$, $W_{\bullet}$, which induces filtrations on the
hypercohomology groups
\[ \H^d(\Omega_{\mZ/\Delta}^{\bullet}(\log D)\otimes \mO_D),\]
which then leads to the limiting mixed Hodge structure on $H_s$. 
We refer to \cite{PetSte} for  a detailed account.\\

If the family $f:\mathcal{Y} \to S$ is defined over $\Q$, there is also a treasure of arithmetical information associated to the situation. We obtain for each 
rational point of $\P^1 \setminus \Sigma$ a Galois-representation on the $l$-adic cohomology 
$H^3_{\textup{{\'e}t}}({Y}_t\otimes_{\Q}\overline{Q},\Q_l)$, and these together make up an $l$-adic sheaf.\\

We will be mainly interested in the case where
\[ h^{12}=1, \]
so the local system $\H$ on $\P^1\setminus \Sigma$ is a so-called
$(1,1,1,1)$-variations and its representations  lands in
$Sp_4(\Z)$. In particular, for each of the $63$ families of {\sc
  Meyer}, we obtain a family of double octics
\[  0=u^2-f_8(x,y,z,w;t).\]
By crepant resolution of the general fibre a family (dropping the
earlier $\hat{}\;$) we obtain such families
\[ f: \mathcal{Y} \to S:=\P^1 \setminus \Sigma \]
and from it an associated $(1,1,1,1)$-variation over $\P^1\setminus \Sigma$.\\

Strictly speaking, an arrangement defined by $f_8=0$ does not define a single family.
If we multiply $f_8$ with a $t$-dependent function $\varphi(t)$, the family defined by
\[0=u^2-\varphi(t)f_8(x,y,z,w;t)\]
is said to be a {\em twist} of the family 
\[0=u^2-f_8(x,y,z,w;t)\] 
Although twisting can be crucially important, its effect is usually easy to analyse, 
and we will consider families differing by a twist as essentially the same.\\

\subsection{Degenerations of (1,1,1,1)-variations}

There are four possibilities for the mixed Hodge diamond of the limiting mixed Hodge structures 
appearing for $(1,1,1,1)$-VHS. The $k$-th row (counted from the bottom) of the diamond gives the 
Hodge numbers of $Gr^W_k$; the monodromy logarithm operator $N$ acts in the vertical direction,  
shifting downwards by two rows. The definition of the weight filtration makes the diagram 
symmetric with respect to reflection in the central horizontal line, whereas complex conjugation 
is a symmetry of the Hodge-diamond along the central vertical axis.
The numbers in each slope $=1$ (so SW-NE-direction) row of the diagram have to add up to the
corresponding Hodge number of the variation, so are all equal to $1$ in our case. 
The cases that arise are:\\

{\bf \em F-point}
\[
\begin{array}{ccccccc}
&&&0&&&\\
&&0&&0&&\\
&0&&0&&0&\\
1&&1&&1&&1\\
&0&&0&&0&\\
&&0&&0&&\\
&&&0&&&\\
\end{array}
\]
In this case $N=0$, so this case occurs if and only if the monodromy is of {\em finite order}. The
limiting mixed Hodge structure is in fact pure of weight three. Often an automorphism of finite order will split the Hodge structure $Gr^W_3$:\\
\[ (1\;\; 1\;\; 1\;\; 1) \mapsto (1\;\; 0\;\; 0\;\; 1)\; +( 0\;\;1\;\;1\;\;0)\]

On the arithmetic side, one expect that when this happens over $\Q$, the characteristic polynomial 
of Frobenius will factor as
\[ (1-a_pT+pT^2)(1-c_pT+p^3T^2)\]
where the $a_p$ and $c_p$ are Fourier coefficients of resp. a weight $2$ and a weight $4$ cusp
form for some congruence subgroup $\Gamma_0(N)$ of $Sl_2(\Z)$.\\
There are also cases where no splitting occurs, but the Euler-factors are determined by a
Hilbert modular form of weight $(4,2)$ for some real quadratic extension of $\Q$.\\  

{\bf \em C-point}
\[
\begin{array}{ccccccc}
&&&0&&&\\
&&0&&0&&\\
&0&&1&&0&\\
1&&0&&0&&1\\
&0&&1&&0&\\
&&0&&0&&\\
&&&0&&&\\
\end{array}
\]
In this case $N \neq 0$, $N^2=0$ and there is a single Jordan block.
The pure part $Gr^W_3$ is a rigid Hodge structure with Hodge numbers
$(1,0,0,1)$. Furthermore, $Gr^W_4$ and $Gr^W_2$ are one-dimensional and
are identified via $N$. This type appears when a Calabi--Yau threefold
acquires one or more ordinary double points, nowadays often called
{\em conifold points}, which explains our name $C$-type point for it.
However, one should be aware that $C$-point do occur not only where
ordinary nodes appear, but also for many other kinds of singularities.\\

On the arithmetical side, one expects to get a $2$-dimensional Galois 
representation $Gr^W_3$ with characteristic polynomial of Frobenius of
the form
\[ 1-a_p T+p^3T^2,\]
where the Frobenius traces are Fourier coefficient of a weight $4$ cusp
form for a congruence sub-group $\Gamma_0(N)$ of $Sl_2(\Z)$, for some level 
$N$.\\

{\bf \em K-point}
\[
\begin{array}{ccccccc}
&&&0&&&\\
&&0&&0&&\\
&1&&0&&1&\\
0&&0&&0&&0\\
&1&&0&&1&\\
&&0&&0&&\\
&&&0&&&\\
\end{array}
\]

In this case we also have $N \neq 0$, $N^2=0$ but there are two
Jordan blocks. In this case the pure part $Gr^W_3 =0$ and $Gr^W_4$,
$Gr^W_2$ are Hodge structures with Hodge numbers $(1,0,1)$, which
are identified via $N$. The Hodge structure looks like that of the
transcendental part of a K3-surface with maximal Picard number, which
explains our name $K$-point for it.\\
 
On the arithmetical side, one expects to get a $2$-dimensional Galois 
representation $Gr^W_2$ with characteristic polynomial of Frobenius of
the form
\[ 1-a_p T+p^2T^2,\]
where the Frobenius traces are Fourier coefficient of a weight 3 cusp
form for a congruence sub-group $\Gamma_0(N)$ of $Sl_2(\Z)$, for some level 
$N$ and character. Such forms always have complex multiplication (CM).
For a nice overview see \cite{Schu1} and \cite{Schu2}.\\

{\bf \em MUM-point}
\[
\begin{array}{ccccccc}
&&&1&&&\\
&&0&&0&&\\
&0&&1&&0&\\
0&&0&&0&&0\\
&0&&1&&0&\\
&&0&&0&&\\
&&&1&&&\\
\end{array}
\]

Here $N^3 \neq 0$ and there is a single Jordan block of maximal size.
The Hodge structures $Gr^W_{2k}$ ($k=0,1,2,3$) are one-dimensional and
necessarily of Tate type. This happens for the quintic mirror at $t=0$
and is one of the main defining properties of Calabi--Yau operators.\\
So at a MUM-point, the resulting mixed Hodge structure is an iterated
extension of Tate--Hodge structures. {\sc Deligne} \cite{Del} has shown
that the instanton numbers $n_1, n_2, n_3, \ldots$ can be seen to 
encode precisely certain {\em extension data} attached to the variation 
of Hodge structures near the MUM-point.\\

\section{Picard-Fuchs operators}

\subsection{Generalities}

By a {\em choice of volume forms on the fibres of} $\overline{f}: \overline{\mY} \to \P^1$ 
we mean a rational section $\omega$ of the relative dualising sheaf
\[ \overline{f}_* \omega_{\overline{\mathcal{Y}}/\P^1} .\]
It restricts to a holomorphic $3$-form $\omega(t)$ on each regular fibre $Y_t$ outside
the divisor of poles and zero's of $\omega$.
If $\gamma(t) \in H_3(Y_t,\Z)$ is a family of horizontal cycles, defined in a 
contractible neighborhood $U$ of $t \in S$, the function
\[ \Phi_\gamma: U \to \C,t \mapsto \int_{\gamma(t)} \omega(t)\] 
is called a {\em period integral} of $f:\mY \to \P^1$. 
It follows from the finiteness of the de\-Rham cohomology group $H^3(Y_t,\C)$ by {\em 
differentiation under the integral sign} that all period functions
$\Phi_{\gamma}(t)$ satisfies the same linear ordinary differential equation, 
called the {\em Picard-Fuchs equation}. The corresponding differential operator 
is called 
{\em Picard-Fuchs operator} and we will write
\[\mP=\mP(\overline{\mY},\omega) \in \Q\langle t,\frac{d}{dt}\rangle .  \]
The order of the operator is clearly at most $\dim H^3(Y_t)$.\\

In the case of double octics  given by an affine equation
\[ u^2-f_8(x,y,z,t) \]
we will always take the volume form 
\[ \omega:=\frac{dxdydz}{u}\]
as three-form and thus the period integrals we are dealing with are
\[ \Phi_{\gamma}(t):=\int_{\gamma(t)} \frac{dxdydz}{u}=\int_{\gamma(t)} \frac{dxdydz}{\sqrt{f_8}}\]
over cycles $\gamma(t) \in H_3(Y_t, \Z)$.

\subsection{Determination of Picard-Fuchs operators}

We have been using two fundamentally different methods to find Picard-Fuchs operators
for concrete examples.\\

{\em Conifold expansion method.}
If in a family of varieties we can locate a vanishing cycle, then the
power series expansion of the period integral can always be computed
algebraically, \cite{CvS1}. The operator is then found from the recursion of the
coefficients. 

Especially for the case of double octics, in many of the $63$ families one can identify  
a {\em vanishing tetrahedron}: for a special value of the parameter one of 
the eight planes passes through
a triple point of intersection, defined by three other planes. In appropriate 
coordinates we can write our affine equation as
\[ u^2=xyz(t-x-y-z)P(x,y,z,t)\]   
where $P$ is the product of the other five planes and depend on a the parameter $t$.
We assume $P(0,0,0,0) \neq 0$.
One can now identify a nice cycle $\gamma(t)$ in the double octic, which consists of 
two parts $\gamma(t)_+,\;(u \ge 0)$ and $\gamma(t)_{-},\;(u\le 0)$ which project onto the 
real tetrahedron $T_t$ bounded by the plane $x=0$, $y=0$, $z=0$, $x+y+z=t$. $\gamma(t)_{+}$
and $\gamma_{-}$ are glued together at the boundary of $T_t$, thus making up a three sphere
in the double octic.  For $t=0$ the tetrahedron and thus the sphere $\Gamma(t)$ shrink to a point. We can write
\[ \Phi_{\gamma}(t)=\int_{\gamma(t)} \omega=2F(t)\]
where
\[F(t)=\int_{T_t} \frac{dxdydz}{\sqrt{(xyz(t-x-y-z)P(x,y,z,t)}}=t \int_{T} \frac{dxdydz}{\sqrt{(xyz(1-x-y-z)P(tx,ty,tz,t)}}, \]
where we used the substitution $(x,y,z) \mapsto (tx,ty,tz)$.
By expanding the integrand in a series and perform termwise integration over the simplex $T:=T_1$,  the
period expands in a series of the form
\[ \Phi_{\gamma}(t) =\pi^2 t (A_0+A_1t+A_2t^2+\ldots)\]
with the coefficients $A_i \in \Q$ if $P_t(x,y,z) \in \Q[x,y,z,t]$ and can be computed explicitly, see \cite{CvS1}. By computing  sufficiently many terms in the expansion, one may 
find the Picard-Fuchs operator by looking for the recursion on the coefficients $A_i$.\\

{\em Cohomology method.} There are many variants for this method, but let us take for sake of
simplicity the case of a smooth hypersurface  $X \subset \P^{n}$. The middle dimensional
 (primitive) cohomology of $X$ can be identified with $H^n(\P^n\setminus X)$ and it elements can
be represented by residues of $n$-forms on $\P^n$ with poles along $X$ of the form
 \[  \frac{P\Omega}{F^k} \]
where $F$ is the defining polynomial for $X$, $\Omega=\iota_E(dVol)$ the fundamental
form for $\P^n$, and $P$ is a polynomial such that the above expression is homogeneous 
of degree $0$. If $P$ runs over the appropriate graded pieces of the Jacobi-ring of $F$, one
obtains a basis 
\[\omega_1,\omega_2,\ldots,\omega_N\]
for $H^n(\P^n\setminus X)$. 
If $F$ depends on an additional parameter $t$, one can differentiate these basis-elements
with respect to $t$ and express the result in the given basis. One obtains thus a
differential system
\[ \frac{d}{dt}
\left (\begin{array}{c} 
\omega_1\\
\omega_2\\
\ldots\\
\omega_N\\
\end{array} \right)=A(t)\left(\begin{array}{c} 
\omega_1\\
\omega_2\\
\ldots\\
\omega_N\\
\end{array} 
\right)
\]
from which one can obtain Picard-Fuchs equations for each $\omega_i$.
This so-called {\em Griffiths-Dwork} method (\cite{Gri1}, \cite{Dw1}) depends 
on the assumption that the hypersurface is smooth, so that the partial derivatives 
$\partial_iF$ form a regular sequence in the polynomial ring. If $X$ has singularities, 
one no longer obtains a basis for the cohomology. Rather one has to determine
the Koszul-homology between the partial derivatives and enter into a spectral
sequence and things become more complicated. In \cite{CvS1} we wrote:\\

 {\em Due to the singularities of $f_8$, 
a Griffiths-Dwork approach is cumbersome, if not impossible.}\\

It was {\sc P. Lairez} who proved us very wrong in this respect.  
His computer program, described in \cite{Lai}, does not aim at finding a 
complete cohomology space, rather it looks for  the smallest space stable under 
differentiation that contains the given rational 
differential form. It does so by going through the spectral 
sequence given by the pole order filtration, where at each step Gr\"obner basis
calculations are done to increase the set of basis forms.\\

For our computations we initially used the method of conifold expansion, 
but  we discovered soon that\\

 {\em Due to the singularities of $f_8$ the conifold expansion approach is cumbersome, if not impossible.}\\

\subsection{Reading the Riemann-symbol}

The amount of information that is contained in the Picard-Fuchs operator 
$\mP$ can not be underestimated. The local system $\H_{\C}$ is isomorphic to local 
system of solutions $Sol(\P)$. Already the local monodromies $T_s$ around the special 
parameter values $s \in \Sigma$ are hard to obtain from topology or a semi-stable 
reduction. But this information can easily be read off from the operator, by studying 
the local solutions in series of the form
\[ t^{\alpha} \sum_{k=0}^N \sum_{n=0}^{\infty}A_{n,k}t^{n}\log(t)^k .\]

There is a delicate interaction between $\mP$ and the Frobenius-operator (see \cite{Dw2}),
so that arithmetical properties of the varieties are tightly linked to $\mP$.
It appears that the Picard-Fuchs operator just abstracts away sufficiently many 
details of the geometry and retains just the right amount of motivic information.\\

We recall that the {\em Riemann-symbol} of a differential operator $\mP \in \C\langle t,\frac{d}{dt}\rangle$
is a table recording for each singular point of the differential operator the
corresponding {\em exponents}, i.e. solutions to the {\em indicial equation} \cite{Ince}.
(In order to have a non-zero series solution of the type described above, one needs that 
$\alpha$ is an exponent at $0$.)\\
  
We found it convenient to express the operators in terms of the logarithmic differentiation
\[ \Theta:=t \frac{d}{dt}\]
and write the operators in {\em $\Theta$-form}
\[ \mP:=P_0(\Theta)+tP_1(\Theta)+t^2P_2(\Theta)+\ldots t^rP_r(\Theta), \;\;\;P_r \neq 0\]
where the $P_i$ are polynomials, in our case of degree four.
The exponents at $0$ are then just the roots of $P_0$, those of $\infty$ the 
roots of $P_r$, with a minus sign. To determine the exponents at other points, 
one just translate this point to the origin, and re-express the operator in $\Theta$-form.\\

The exponents capture the semi-simple part of the monodromy at the corresponding singular point. 
The logarithmic terms appearing in the solutions encode the Jordan structure of the unipotent part. 
In general logarithmic terms  may appear between solutions with integer
difference in exponents, but in the geometrical context, as a rule, a logarithm appears {\em always} 
 precisely when two exponents become {\em equal}.\\

A {\bf C-point} can be expected when the exponent are of the form
\[ \alpha-\epsilon\;\;\; \alpha\;\;\; \alpha\;\;\; \alpha+\epsilon \]
The archetypical case is $0\;\; 1\;\; 1\;\; 2$, indicating the presence of local solutions of the form
\[ \phi_0=1+a_1t+\ldots, \phi_1(t)=t+b_1t+\ldots, \phi_2(t)=log(t)\phi_1(t)+c_1 t+\ldots, \phi_3=t^3+d_1t^4+\ldots\]
as appear in a pure conifold smoothing.\\ 

A {\bf $K$-point} can be expected to appear when the exponents are of the form
\[ \alpha\;\;\; \alpha\;\;\; \beta\;\;\; \beta,\;\;\;(\alpha <\beta) .\]
The archetypical case is $0\;\; 0\;\; 1\;\; 1$, indicating the presence 
of local  solutions of the form
\[ \phi_0=1+a t+\ldots, \phi_1(t)=\log(t)\phi_0(t)+b t+\ldots, 
\phi_2(t)=t+c t^2+\ldots, \phi_3=\log(t)\phi_1(t)+d t^2+\ldots\]

A {\bf MUM-point} can be expected to appear, when all exponents are equal
\[ \alpha\;\;\; \alpha\;\;\; \alpha\;\;\; \alpha\;\;\;\]
The archetypical case is $0\;\;0\;\;0\;\;0$, with the famous Frobenius
basis of solutions of the form
\[\phi_0=1+at+\ldots,\phi_1(t)=\log(t)\phi_0(t)+\ldots,\phi_2(t)=\log(t)^2\phi_0(t)+\ldots, \Phi_3(t)=\log^3(t)\phi_0(t)+\ldots\]

An {\bf $F$-point} can be expected in all other cases. If all exponents
are integral (and no logarithms appear) we have trivial monodromy and we speak
of an {\em apparent singularity}, if the exponents are {\em not} 
$0\;\;1\;\;2\;\;3$, which would be the exponents at a regular point. 
The archetypical apparent singularity is signaled by the exponents
\[ 0 \;\;\;1 \;\;\;3\;\;4 .\] 
In general, we will call any $F$-point {\em with non-equally spaced exponents} 
an $A$-point.\\ 
When we write the operator $\mP$ in the form
\[\frac{d^4}{dt^4} +a_1(t) \frac{d^3}{dt^3}+a_2(t) \frac{d^2}{dt^2}+a_3(t) \frac{d}{dt}+a_4(t)
\]
where $a_i(t) \in \C(t)$, then the function
\[ Y(t):=e^{-\frac{1}{2} \int a_1(t) dt}\]
is called the {\em Yukawa coupling} and its (simple) zero's typically are apparent singularities with exponents $0\;\; 1\;\; 3\;\; 4$.\\
 
{\bf \em Basic transformation theory}

In what follows, we will not distinguish between an operator $\mP$ and 
the operator $\mP' =\varphi(t) \mP$ obtained from multiplying $\mP$ by a 
rational function $\varphi(t)$, as they determine the same local system of 
solutions on $\P^1 \setminus \Sigma$. (Of course, in the finer theory of 
$\mathcal{D}$-modules one has to  distinguish very well between $\mP$ and $\mP'$).\\ 

Often one has to make simple transformations on differential operators.\\

1) The simplest are those induced by M\"obius transformations, coming from
fractional linear transformations
\[ t \mapsto \frac{a t+b}{ct+d}\] 
of the coordinate in $\P^1$. Of course, this just changes the position of
the singular points, the corresponding exponents are preserved. We call 
two operators related in this way {\em similar operators}.\\

2) If $\omega(t)$ and $\omega'(t)$ are two different choices of volume
form on the fibres, then 
\[\omega'(t) =\varphi(t)\omega(t) \]
where $\varphi(t)$ is a rational function. The corresponding Picard-Fuchs 
operators $\mP(\mY,\omega)$ and $\mP(\mY,\omega')$ will be related in a
certain way.\\
 
More generally, if $y(t)$ satisfies $\mP y(t)=0$, and $\phi(t)$ is a 
rational function, then the function $Y(t):=\varphi(t)y(t)$ 
will satisfy another differential equation $\mQ Y(t)=0$ that is rather
easy to determine.  We will say that $\mQ$ is {\em strictly equivalent} 
to $\mP$. Its effect on the Riemann symbol will be a shift of all exponents 
by an amount given be the order of $\varphi(t)$ at the point in question.
For example, the effect of multiplication by $t$ shifts the exponents
at $0$ one up, those at $\infty$ one down.\\

3) As already mentioned above, if $\phi(t)$ is a rational function of $t$, 
then the families of double octics
\[ u^2=f_8,\;\;\;\textup{and}\;\;\; u^2=\phi(t)f_8\]
are said to differ by a {\em twist}.
Replacing $\phi(t)$ by $\phi(t)\varphi(t)^2$ does not change the fibration 
birationally, as can be seen by replacing $u$ by $\varphi(t)u$. 
The volume form 
\[ \omega:=\frac{dxdydz}{\sqrt{f_8}} \]
for $u^2-f_8$ and 
\[ \omega':=\frac{dxdydz}{\sqrt{\phi(t)f_8}} \]
for  $u^2-\phi(t) f_8$ differ by the square root of a rational function
\[ \omega=\sqrt{\phi(t)}\omega'\]  

More generally, if  $y(t)$ satisfies a differential equation $\mP y(t)=0$, then 
$Y(t):=\phi(t)y(t)$, where $\phi(t)$ is an algebraic function, will satisfy
another differential equation $\mQ$ that is easy to determine, knowing only
$\mP$. We will then call $\mP$ and $\mQ$ {\em equivalent}.\\

Its effect on the Riemann symbol is also rather easy to understand.
If $\phi(t)$ has near $a$ the character of $(t-a)^{\epsilon}$, then the exponents 
at $a$ get all shifted by the amount $\epsilon$:

\[ \alpha,\beta,\gamma,\delta \mapsto \alpha+\epsilon,\beta+\epsilon,\gamma+\epsilon,\delta+\epsilon\]

4)  If we replace $t$ by $t=\psi(s)$ for some function $\psi(s)$ we can rewrite the operator the
operator $\mP$ in terms of the variable $s$ and obtain an operator $\psi^*\mP$ in $s, \frac{d}{ds}$
that we call the {\em pull-back} of $\mP$ along the  map $\psi$. Very common are pull-backs
by the map $t=s^n$, which geometrically is an $n$-fold covering map of $\P^1$, with total
ramification at $0$ and $\infty$. This operation leads to a division of the exponents at $0$ and $\infty$:
\[ \alpha,\beta,\gamma,\delta \mapsto \alpha/n,\beta/n,\gamma/n,\delta/n\]

The most general transformations one has to allow are those given by algebraic coordinate 
transformations, which are multi-valued maps from $\P^1$ to itself, which properly understood 
are given by {\em correspondences} via a smooth curve $C$:
\[
\begin{diagram}
\node{C} \arrow{e,t}{p}\arrow{s,l}{q} \node{\P^1}\\
\node{\P^1}
\end{diagram}
\] 
and $\psi^*$ 'is' $q_*p^*$, which just means that
\[ p^*\mP = q^*\mQ\]
In such a case  we will say that $\mP$ and $\mQ$ are {\em related operators}.
The effect of these transformations on the Riemann-symbol can be traced back 
the local ramification behaviour of $p$ and $q$; we will not spell out 
the details.\\ 

It is easy to see that under pull-back one can not get rid of a $MUM$, $K$, $C$ or
$A$-point. Only an $F$-point with equidistant exponents may turn into the 
non-singularity with exponents $0,1,2,3$. 
Note in particular that an operator with a MUM-point can not be related to an
operator without a MUM-points, etc.

\section{Orphans of order $2$}

It turns out that there are seven arrangements that lead to a second order operator. These are the arrangements
\[{\bf 4},\;\; {\bf 13},\;\; {\bf 34},\;\; {\bf 72},\;\; {\bf 261},\;\; {\bf 264},\;\; {\bf 270} .\]
The differential equations for each of these cases was computed; the results are recorded in Appendix B. All operators turn out to be of a very simple type, directly related to the Legendre differential equation, which is the hypergeometric equation
\[\Theta^2-16 t(\Theta+1/2)^2\]
and Riemann symbol
\[\left\{
\begin{array}{ccc}
0&1/16&\infty\\
\hline
0&0&1/2\\
0&0&1/2\\
\end{array}
\right\}.
\]
This also is the Picard-Fuchs operator of the elliptic surface with 
Kodaira fibres $I_2,I_2,I_2^*$.
The differential equations in the cases ${\bf 72}$ and ${\bf 270}$ are a 
bit different.

At first sight it is very surprising to find a second order equation for
such octic triple integrals. As explained in \cite{Roh1}, the appearance
of a certain {\em maximal automorphism} causes the Picard-Fuchs operator 
to be of order two. In \cite{CKC} such maximal automorphism were identified
in five of the seven cases. For the remaining two cases {\bf 264} and 
{\bf 270} we do not have such a simple explanation for the appearence of 
a second order Picard-Fuchs equation.\\

A priori, there seem to be two different scenario's  in which the
Picard-Fuchs operator of a family can reduce to an operator of order two. 
It could happen that the $(1,1,1,1)$-VHS splits as sum into a (rigid) 
$(1,0,0,1)$ Hodge structure and a variable $(0,1,1,0)$, coming from an 
elliptic curve. Or it could be that the $(1,1,1,1)$-VHS is a tensor product 
of a constant $(1,0,1)$-Hodge structure with a variable $(1,1)$-VHS coming 
from a family of elliptic curves. This $1,0,1$ should be the transcendental 
part  of $H^2$ of a K3-surface with Picard number $20$.\\
In our situation it is always the second alternative that has to occur, 
as we are looking at the Picard-Fuchs operator for the period integrals of
the holomorphic volume form $\omega$, which in the first case would be
constant.\\

It is of interest to identify the K3-surface in the geometry of the
arrangement. For example, for the octic corresponding to arrangement {\bf 13} 
\[Y_t: 0=u^2 - x y z (x+y) (y+z) w (x-z-w) (x-z-tw)\]
one can understand its relation to the K3-surface 
in the following way.
Replacing $u$ by $u/(x-z)$ we see the double octic is equal to
the normalisation of the {\em double dectic}
\[u^2=xyz(x+y)(y+z)(x-z)w(x-z-w)(x-z-tw)(x-z) .\]
The first $6$ factors now only depend on the variables $x,y,z$
and determine a double sextic K3-surface $S$ with equation
\[ p^2=xyz(x+y)(y+z)(x-z)\] 

It is the famous {\em most algebraic K3-surface} \cite{Vin}, which comes
with the weight $3$ modular form named $16$ in appendix A.

\begin{center}
\includegraphics[height=5cm]{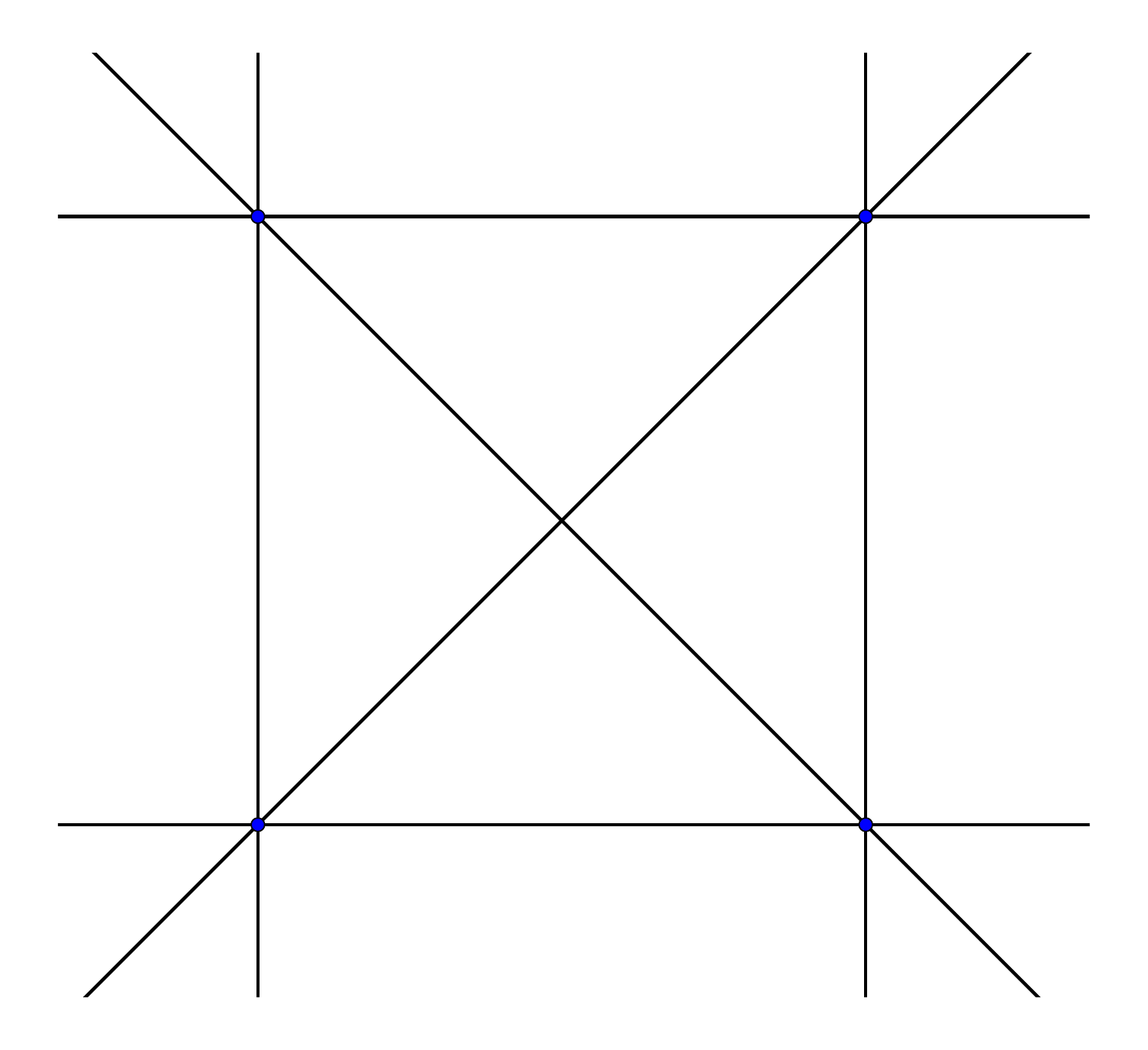}
\end{center}

The last four factors only depend on $w$ and $\xi:=x-z$, with $t$ as
parameter and determine a double quartic family of elliptic curves $E_t$
given by the equation
\[ q^2=w(\xi-w) (\xi-tw) \xi\]
By dividing out the involution $\iota$ induced by  $p\mapsto -p$, $q \to -q$ 
acting on $S \times E_t$ we get back our double dectic, $u=pq$.
Hence we see that our original double octic is birational to
\[ S \times E_t/\iota \lra Y_t \]
so that we see that $Y_t$ is a simple instance of the {\em Borcea-Voisin construction}, see \cite{CM}.\\

It turns out that, unexpectedly, in all cases except ${\bf 270}$ we end up with 
the modular form $16$ as constant factor.

In case ${\bf 270}$ we get modular form $8$, attached to the
double sextic K3-surface 
\begin{center}
\includegraphics[height=5cm]{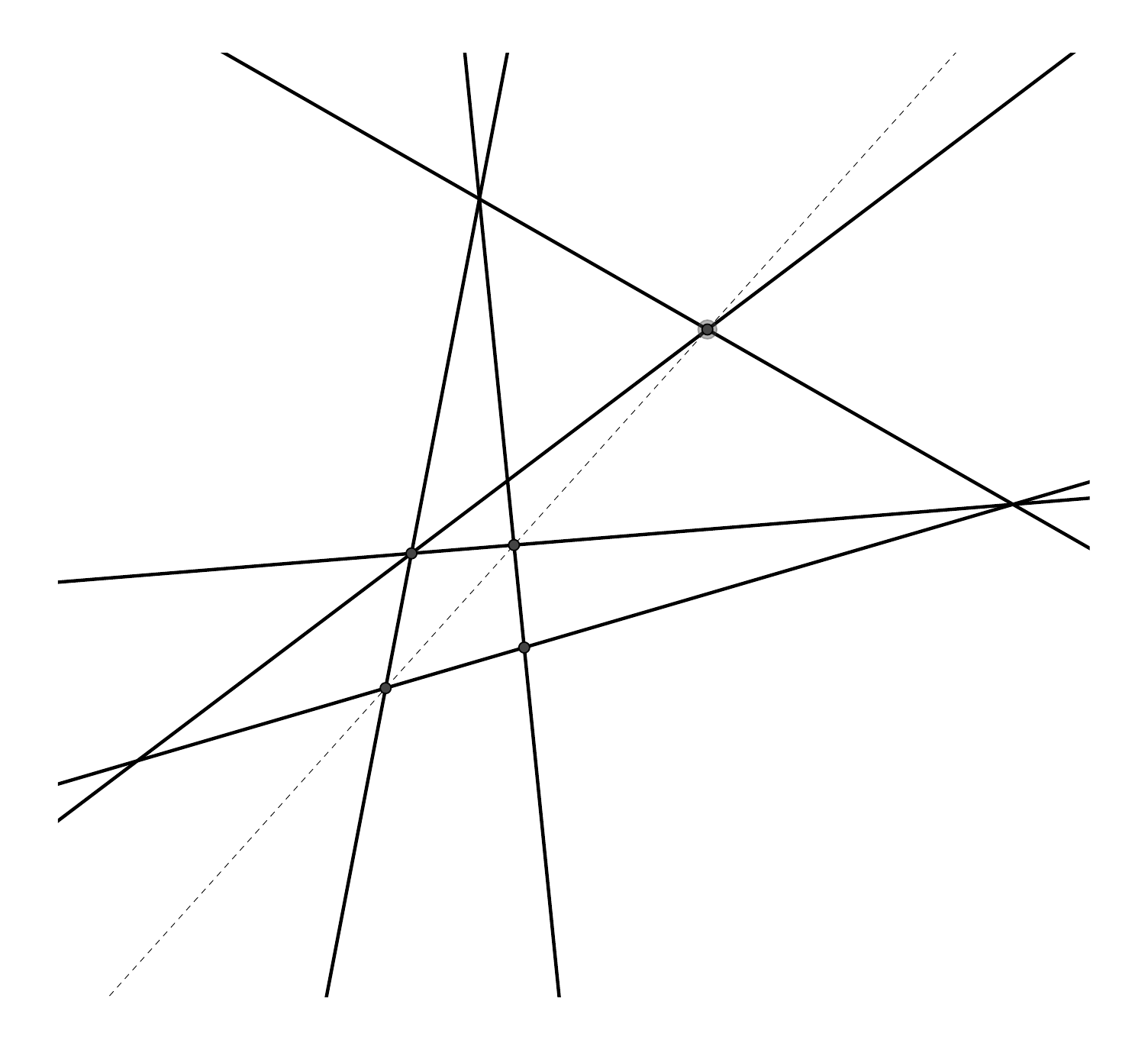}
\end{center}

These two K3-surfaces appear in a nice pencil; in \cite{GT} one finds a
very detailed account of their geometry and arithmetic.\\

\section{Orphans of order $4$}
Of the $63$ families it tuns out that $18$ are fourth order orphans. These
are much more interesting and are collected in Appendix C. We sort them  
according to types of singularities that appear.

\subsection{The two KKCC-operators}

It turns out that there are two different but similar operators with 
two points of type $K$ and two points of type $C$. In each case there is
a pair of arrangements related to them.\\
 
{\bf \em The arrangements ${\bf 33}$ and ${\bf 70}$ with $h^{11}=49$.}\\ 

These two arrangements are birational, via the map
\[(x,y,z,v)\longmapsto (tyv, xy-yz, tzv-xz+{z}^{2}, {x}^{2}-xz-txv).\]

(Here and below we will always refer to the equations of the arrangements
given in Appendix C.)

The Picard-Fuchs operator for ${\bf 33}$ has Riemann symbol
\[
\left\{
\begin{array}{cccc}
0&1&2&\infty\\
\hline
0&0&0&1/2\\
0&1/2&1&1/2\\
1&1/2&1&3/2\\
1&2&2&3/2\\
\end{array}
\right\}
\]
from which we see that $0$ and $\infty$ are $K$-points, and $1$ and $2$
are $C$-points. The operator for $\bf{70}$ differs from it by $t \mapsto -t$. 
This operator was obtained in \cite{CvS1} by conifold expansion.

As explained above, to each rational $K$-point there is attached a weight three modular 
form and to each rational $C$-form a weight four modular form. It is well-known that these 
forms can be determined by counting points over finite fields and we will omit all details on their calculation. 

We found it convenient to write the names of these forms {\em above} the corresponding points 
of the Riemann-symbol, as to form a {\em decorated Riemann-symbol.}

\[
\left\{
\begin{array}{cccc}
8&32/2&8/1&16\\
\hline
0&1&2&\infty\\
\hline
0&0&0&1/2\\
0&1/2&1&1/2\\
1&1/2&1&3/2\\
1&2&2&3/2\\
\end{array}
\right\}
\]
(For the naming of the modular forms we are using the reader is referred to Appendix A.)
We see that two different types of $K3$-surfaces, corresponding to the forms $8$ and $16$ should 
appear in the semi-stable model of the singular fibres. This is in accordance with the fact that 
there is no symmetry that fixes the $C$-points and interchanges the $K$-points.\\

{\bf \em The arrangements ${\bf 97}$ and ${\bf 98}$ with
  $h^{11}=45$.}\\

This is the other pair of arrangements that also lead  to an operator with two $K$ and two $C$-points. The arrangements ${\bf 97}$ and ${\bf 98}$ are also birational, via the map\\

\[(x,y,z,v)\longmapsto\left( (tv-z-v)  \left( x+y+z+v \right) ,txz,( -tv+z+v) y, ( tv-z-v)  ( x+y) \right) \]

The decorated Riemann symbol of the operator for ${\bf 97}$ is

\[
\left\{
\begin{array}{cccc}
32/1&8&8/1&8\\
\hline
0&-1&-2&\infty\\
\hline
  0&0&0&1/2\\
1/2&0&1&1/2\\
1/2&1&1&3/2\\
  1&1&2&3/2\\
\end{array}
\right\}
\]
By shifting the exponents by $1/2$ in $0$ and $-1/2$ in $\infty$ followed by a simple M\"obius
transformation one can transform it to operator ${\bf 98}$.\\
 
We again see two $K$ points, but this time we find that at both of them the modular form is $8$, 
which suggests that the operator of {\bf 98} has a symmetry that interchanges the two $K$-points. 
Indeed, by shifting the finite $C$-point to the origin,
one finds that the operator is symmetric under $t \mapsto -t$ and thus can 
be pulled-back by the squaring map from the nice operator
\[ \mathcal{A}:=\Theta^{2}( \Theta-1)^{2}+t\, \Theta^{2}(32\,\Theta^{2}+3)+4\,t^2 \left( 4\,\Theta+1 \right)  \left( 2\,\Theta+1 \right) ^{2} \left( 4\,\Theta+3 \right) \]
which has extended Riemann symbol
\[
\left\{
\begin{array}{cccc}
8&8/1&32/1\\
\hline
0&-1/16&\infty\\
\hline
0&0&1/4\\
0&1/2&1/2\\
1&1/2&1/2\\
1&1&3/4\\
\end{array}
\right\}
\]

The symmetry of the operator is also visible as a symmetry of the arrangements.\\

\subsection{The two KCCC operators}

It turns out that there are also two different operators with a single $K$ point, 
but with three additional $C$-points. The first of these operators is related to two essentially different pairs of double octic arrangements, namely\\

{\bf \em The arrangements $\bf 35$ and  $\bf 71$ with $h^{11}=49$.}\\
and\\
{\bf \em The arrangements $\bf 247$ and $\bf 252$ with $h^{11}=37$.}\\

As the Riemann-symbol suggests, the Picard-Fuchs operators for these four cases are 
related by a simple M\"obius transformation and multiplication with an algebraic function
and we will only analyse one of the cases. 

We know that ${\bf 247}$ and ${\bf 252}$ are birational arrangements, but we were unable
to find a birational transformation between ${\bf 35}$ and ${\bf 71}$. The coincidence of the Picard-Fuchs operators (up to transformation) strongly suggest that there exist a
{\em correspondence} between ${\bf 35}$ and ${\bf 247}$, but again we were unable to find
it. This is an illustration of the power of Picard-Fuchs operators to make geometrical
predictions.\\   

The decorated Riemann symbol of ${\bf 35}$ is
\[
\left\{
\begin{array}{cccc}
8/1&8/1&8&8/1\\
\hline
-1&0&1&\infty\\
\hline
0&  0&0&1/2\\
1&1/2&0&1\\
1&1/2&1&1\\
2&  1&1&3/2\\
\end{array}
\right\}
\]
so at all $C$-points we find the modular form $8/1$.

Indeed, the operator of $\bf 35$ has symmetries interchanging the $C$-points.
If we shift the exponents at $-1$ by $1/2$, at $\infty$ by $-1/2$, and then
bringing  the $K$-point to $0$ and the C-point at $-1$ to $\infty$, we obtain 
an operator symmetric under $t \mapsto -t$, so it is seen to be pull-back
by the squaring map of 
\[ \mathcal{B}:=
{\Theta}^{2} \left( 2\,\Theta-1 \right) ^{2}+t \left( 4\,{\Theta}^{2}+2\,\Theta+1 \right) 
 \left( 4\,\Theta+1 \right) ^{2}+ t^2 \left( 4\,\Theta+1 \right)  \left( 4\,\Theta+3
 \right) ^{2} \left( 4\,\Theta+5 \right) \]
with decorated Riemann-Symbol
\[
\left\{
\begin{array}{ccc}
8&8/1&8/1\\
\hline
0&-1/8&\infty\\
\hline
0&0&1/4\\
0&1/2&3/4\\
1/2&1/2&3/4\\
1/2&1&5/4
\end{array}
\right\}
\]
At both conifold points we have modular form $8/1$, and it turns out that 
there is a further symmetry in the operator that exchanges these and thus 
can be obtained as pull-back from yet another operator with three singular 
points: 
\[{\Theta}^{2} \left(4\,\Theta - 1 \right) ^{2}+ 2\,t \left(8\,\Theta+1 \right) 
 \left(32\,{\Theta}^{3}+28\,{\Theta}^{2}+19\,\Theta+4 \right) + t^2  \left(8\,\Theta+1
 \right)  \left( 8\,\Theta+9 \right)  \left( 8\,\Theta+5 \right) ^{2} \]
Its decorated Riemann symbol is
\[
\left\{
\begin{array}{ccc}
8&?&8/1\\
\hline
0&-1/8&\infty\\
\hline
0&0&1/8\\
0&1/2&5/8\\
1/4&1&5/8\\
1/4&3/2&9/8\\
\end{array}
\right\}
\]
Note that now there appears an singularity at $-1/8$ with monodromy of
order $2$. So the Calabi-Yau threefold appearing at this point is special,
but the $?$ indicates that we were not able yet to identify any modular 
forms.\\

{\bf \em The arrangements $\bf 152$ and $\bf 198$ with $h^{11}=41$}.\\

These two arrangement give rise to another operator with one $K$ and three $C$-points.
The coincidence of Hodge numbers and Picard-Fuchs operator (up to equivalence) suggest,
that the two arrangements are birational. Again, we were unable to find the map.\\ 

The decorated Riemann symbol of $\bf 152$ is 
\[ 
\left\{
\begin{array}{cccc}
8&8/1&32/1&8/1\\
\hline
-1&0&1&\infty\\
\hline
0&0&0&1/2\\
1/2&1/2&0&1\\
1/2&1/2&2&1\\
1&1&2&3/2\\
\end{array}
\right\}
\]
The operator of ${\bf 198}$ is equivalent to it: by shifting the exponents at $-1$ by $1/2$ and at $\infty$ by $-1/2$. 
 
The appearance of the forms $8/1$ at both $0$ and $\infty$ suggests that there is a symmetry interchanging these points. And indeed, it turns out that the operator for ${\bf 152}$ is also a pull-back from the operator 
$\mathcal{A}$!

\subsection{The ACCC-operator}
\vskip10pt
{\bf \em The arrangements 153, 197 with $h^{11}=41$}\\ 

The arrangements $\bf 153$ and $\bf 197$ are birational, via the map\\
\[ (x,y,z,v)\longmapsto \left((  xt+vt-y ) tv,
  (  xt+vt -y-z ) y, ( x+v ) tz,
  (  xt+vt-y ) tx\right)\]

The decorated Riemann symbol for $\bf 153$ is:
\[
\left\{ 
\begin {array}{cccc}
32/1&32/2&8/1&32/2\\
\hline
0&-1&-2&\infty\\
\hline
0&0&0&1/2\\
1/2&1/2&1/2&1\\
1/2&1/2&3/2&1\\
1&1&2&3/2\\
\end {array} \right\} 
\]
The point $t=-2$ is very remarkable. By expanding the solutions around
the singular point $-2$ one sees that the monodromy is of order two; no
logarithmic terms arise. As a result, the corresponding limiting MHS remains
pure of weight three. On the other hand, the arrangement at $t=-2$
specialises to the rigid arrangement $93$ and the corresponding double 
octic has a rigid Calabi-Yau (with modular form $8/1$) as resolution.
This implies that any semi-stable fibre at $t=-2$ needs to have 
further components to account for the change in cohomology between
special fibre at $-2$ and general fibre. We note that the theorem
of {\sc Kulikov} implies that a similar phenomenon can not happen for 
K3-surfaces. This and other examples will be studied in more detail in
future paper.\\  

As suggested by the modular forms, there could be a symmetry interchanging the two $32/2$ points. By first shifting exponents by $1/2$ at $-2$ and $-1/2$ at $\infty$ and then bringing $-2$ to $\infty$, the operator is pulled back from 
the nice operator with three singular points:

\[\mathcal{C}:=
\Theta(4\Theta-1)(2\Theta-1)-t(4\Theta+1)^2(4\Theta^2+2\Theta+1)+t^2(4\Theta+1)
(2\Theta+1)(4\Theta+5)(\Theta+1)\]

with extended Riemann symbol

\[\left\{ \begin {array}{ccc}
32/1&32/2&8/1\\
\hline
0&1&\infty\\
\hline   
  0&  0&1/4\\ 
1/4&1/2&3/4\\  
1/4&1/2&3/4\\ 
1/2&  1&5/4\\
\end {array} \right\} .
\]

\subsection{The KCCCC-operator}

{\em The arrangement $\bf 243$} ($h^{11}=39$) 
leads to the rather complicated operator
\[
\Theta\, \left( \Theta-2 \right)  \left( \Theta-1 \right) ^{2}
-\frac16\,t\Theta\, \left( \Theta-1 \right)  \left( 19\,{\Theta}^{2}-19\,\Theta+9 \right) 
+\frac13\,{t}^{2}{\Theta}^{2} \left( 11\,{\Theta}^{2}+4 \right)\]
\[ -\frac1{24}\,{t}^{3} \left( 11\,{\Theta}^{2}+11\,\Theta+5 \right)
\left( 2\,\Theta+1 \right) ^{2} 
+\frac1{48}\,{t}^{4} \left( 2\,\Theta+3 \right) ^{2} \left( 2\,\Theta+1 \right) ^{2}
\]

Its decorated Riemann symbol is

\[
\left\{ 
\begin {array}{ccccc}
12/1&32/2&6/1&8/1&8\\
\hline
0&1&\frac32&2&\infty\\
\hline
0&  0&0&0&1/2\\
1&1/2&1&1&1/2\\
1&1/2&1&1&3/2\\
2&1  &2&2&3/2\\
\end {array} \right\} 
\]

We claim it can not be simplified further, as the modular forms at the cusps are all different.

\subsection{The ACCCCK-operator}

{\bf \em Arrangement ${\bf 250}$ and ${\bf 258}$} define
birational families of Calabi-Yau threefolds with $h^{11}=37$.
The decorated Riemann-symbol of ${\bf 250}$ is
\[
\left\{ 
\begin {array}{cccccc}
6/1&8/1&h&8/1&6/1&8\\
\hline
-2 &-1&-1/2&0&1&\infty\\
\hline
0&  0&0&  0&0&1/2\\
1&1/2&1&1/2&1&1/2\\
1&1/2&3&1/2&1&3/2\\
2&1  &4& 1 &2&3/2\\
\end {array} \right\} 
\]
Very remarkably, at the apparent singularity at the point $-1/2$ there appears the Hilbert-modular form $h$ of level 
$6\sqrt{2}$ and weight $(4,2)$. When we shift this apparent singularity to the 
origin, we obtain an operator that is symmetric with respect to the involution 
$t \mapsto -t$. But there is no obvious corresponding symmetry in the family. 
Counting points, it appears that the number of points  of the fibres at  $t$ and $-t$ are 
equal or opposite mod $p$, according to the quadratic character $\left( \frac{2}{p} \right)$. 
In fact the transformation
\[\scriptsize
\left(  \begin{array}{c}
    x\\y\\z\\t
  \end{array}\right)
\longmapsto
\left(\begin{array}{l}  
2\, \left( -y-z+v \right) x \left( x+v+y+z \right)  \left(  \left(
    t+1/2 \right) y-v-x-z \right) \\
4\, \left( x+v+y+z \right)  \left( -1/2\,{z}^{2}+ \left( v/2-x/2-y \right) z-1/2
  \,{y}^{2}+ \left( v/2-x/2 \right) y+vx \right) z\\
\left( {y}^{2}+ \left( -v+x+2z \right) y + {z}^{2}+ \left( x-v \right) z-2vx \right)  \left( x+v+y+z \right) 
\left(  \left( t+1/2 \right) y-v-x+z \right) \\
\left( y+z \right)  \left( v-x-y-z \right) ^{2} \left( ty-v-x+y/2-z \right) 
\end{array}\right)
\]
gives a correspondance between a fiber of the family and the quadratic
twist by 2 of the opposite fiber.
In particular, the modular forms labelled $8/1$ actually occuring differ by this character.
From this state of affairs it seems natural that at the symmetry point the Hilbert modular form for $\sqrt{2}$ appears. We plan study this example more
carefully in a future paper.\\

We remark that although the modular forms at corresponding fibres $0$ and $-1$ are in both cases $8/1$, they correspond to two non-birational rigid Calabi-Yau configuration ${\bf 69}$ ($h^{11}=50$) and ${\bf 93}$ ($h^{11}=46$).\\

Similarly, at $1$ and $-2$ we have modular form $6/1$, but rigid Calabi-Yau configurations ${\bf 245}$ ($h^{11}=38$)
resp. ${\bf 240}$ ($h^{11}=40$). So it seems improbable that there is a birational map relating the fibre at $t$ to the fibre at $-1/2-t$. Geometrically, the symmetry of the Picard-Fuchs operator is surprising.\\

Using the symmetry we see that the operator is pulled back from a simpler
operator with the following Riemann symbol:
\[
\left\{ 
\begin {array}{cccc}
h&8/1&6/1&\infty\\
\hline
0&1&9&\infty\\
\hline
0  &  0&  0& 1/4\\
1/2&1/2&  1& 1/4\\
3/2&1/2&  1& 3/4\\
2  &  1&  2& 3/4\\
\end {array} \right\} 
\]

\subsection{The KCCCCC-operator}

{\bf \em Arrangement ${\bf 248}$ with $h^{12}=37$.}

Decorated Riemann symbol

\[
\left\{
\begin{array}{cccccc}
12/1&6/1&16&6/1&12/1&32/1\\
\hline
-2&-3/2&-1&-1/2&0&\infty\\
\hline
0&0&0&0&0&1/2\\
1&1&0&1&1&1\\
1&1&2&1&1&1\\
2&2&2&2&2&3/2\\
\end{array}
\right\}
\]

There is a symmetry pairing the points with same modular form
and fixing the points belonging to the modular forms $16$ and $32/1$.
First shift the $K$-point to $0$, then the operator is seen the be
pull-back via quadratic map from an operator with Riemann symbol
\[
\left\{
\begin{array}{cccc}
16&12/1&6/1&32/1\\
\hline
0&1&1/4&\infty\\
\hline
0&0&0&1/4\\
0&1&1&1/2\\
1&1&1&1/2\\
1&2&2&3/4\\
\end{array}
\right\} .
\]

This symmetry, after having discovered it from the operator,  can be seen in the arrangement.\\

\subsection{The ACCCCCC-operator}

{\bf \em Arrangements $\bf 266$ and ${\bf 273}$ } have both $h^{11}=37$ and have rather complicated
but {\em identical} Picard-Fuchs equations, with Riemann symbol

\[\left\{
\begin{array}{cccccccc}
6/1&32/1&6/1&32/1&6/1&?&?&32/1\\
\hline
-1&-1/2&-1/4&0&1/2&(-1+\sqrt{-3})/4&(-1-\sqrt{-3})/4&\infty\\
\hline
0&0  &0&0  &0&0&0&1/2\\
1&1/2&1&1/2&1&1&1&1\\
1&1/2&1&1/2&1&3&3&1\\
2&1  &2&1  &2&4&4&3/2\\
\end{array}
\right\}
\]

Arrangement ${\bf 266}$ and ${\bf 273}$ both have $8$ quadruple points and
no five-fold point or triple lines. The difference between the two is rather 
subtle; ${\bf 266}$ contains six planes which are in general position, which 
is not the case for ${\bf 273}$. In ${\bf 266}$ there are six of the quadruple 
points in a plane, which is not the case for ${\bf 273}$, etc. So clearly the 
two configurations are projectively very different. On the other hand, they 
have the same Hodge number $h^{11}=37$, and the equality of their Picard-Fuchs 
operators clearly suggest the varieties belonging to the two arrangements are 
birational, but were unable to find any transformation.\\

The operator can be reduced to one with three singular points via the following
steps
\begin{itemize}
\item Make the exponents at all $32/1$-points equal to , by exponent-shift.
\item Translate the A-points to $0$ and $\infty$. Of course, the
resulting operator has coefficients in $\Q(\sqrt{-3})$.
\item Make a pull-back of order three, i.e. we use $t^3$ as new coordinate.
\item The operator has now four singular points, where the points at 
$0$ and $\infty$ have exponents $0,1/3,1,4/3$. Bring the two other points
to $0$ and $\infty$.
\item The result is an operator invariant under $t \mapsto -t$. Make
a quadratic pull-back and bring the $A$ point to the origin. 
\item The result is, after a scaling of the coordinate, the following operator

\[
4 \Theta (3 \Theta - 1)(\Theta-1) (3 \Theta - 4)
+6 t \Theta (3 \Theta - 1) (288 \Theta^2  - 96 \Theta + 35)
+144 t^2 (12 \Theta +1)(3 \Theta+1)^2 (12\Theta+7) 
\]

with decorated Riemann symbol
\[
\left\{
\begin{array}{ccc}
  ?&6/1&32/1\\
\hline
  0&-1/36&\infty\\
\hline
  0&  0&1/12\\
1/3&1/2&1/3\\
  1&1/2&1/3\\
4/3&  1&7/12\\
\hline
F&C&C\\
\end{array}
\right\}
\]
\end{itemize}

\section{Birational nature of the Picard-Fuchs operator}

\subsection{Birational maps and strict equivalence}

We would like to formulate a statement expressing the idea that the 
Picard-Fuchs operator is in some sense a birational invariant of a family. 
We formulate three theorems to this effect.  

{\bf Theorem 1:} Consider $\mathcal{X} \lra S$ and $\mathcal{X'} \lra S$
two proper smooth families of Calabi-Yau varieties over $S=\P^1 \setminus \Sigma$. Let $\omega$ and $\omega'$ be holomorphic volume forms on $\mathcal{X}$ and
$\mathcal{X'}$. If there is a birational map
\[ \phi:\mathcal{X} \lra \mathcal{X'},\]
then the Picard-Fuchs operators $\mP(\mathcal{X},\omega)$ and $\mP(\mathcal{X'},\omega')$ are strictly equivalent.\\

{\bf proof:}\\
Let $X$ and $X'$ be fibres over the same general point.
Recall that V. Batyrev \cite{Bat2} proved that if
\[\phi:X \lra X'\]
is a birational map between smooth Calabi-Yau manifolds, then it induces 
isomorphisms
\[ \phi^*: H^{p,q}(X') \to H^{p,q}(X)\]
of Hodge groups. As a consequence, $\phi$ induces isomorphisms
\[\phi_*: H_n(X,\Q) \lra H_n(X',\Q), \;\;\;\phi^*: H^n(X',\Q) \lra H^n(X,\Z) \]
If $\omega \in H^{n,0}(X)=H^0(X,\Omega_X^n)$ and  $\omega' \in H^{n,0}(X')=H^0(X',\Omega_X'^n)$ are holomorphic volume forms on $X$ and $X'$, then
\[ \phi^*(\omega') =\phi \omega\]
for some $\phi \in \C^*$. Furthermore, one has
\[\int_{\gamma} \phi^*(\omega') =\int_{\phi_*(\gamma)} \omega \]

In the relative situation we consider a birational map
that induce fibrewise birational isomorphism
\[ \phi: \mathcal{X}  \dashrightarrow \mathcal{X'} \]
If $\omega(t)$ and $\omega'(t)$ are volume forms on $\mathcal{X}$ resp. 
$\mathcal{X'}$, then
\[ \phi^*(\omega'(t)) =\varphi(t) \omega(t)\]
for some rational function $\varphi(t)$
So we have
\[\varphi(t)\int_{\gamma(t)}\omega(t)=\int_{\gamma(t)} \phi^*(\omega'(t)) = \int_{\phi_*(\gamma(t))} \omega'(t) \]
which shows that the period integrals for $\mathcal{X}$ and $\mathcal{X'}$ differ by multiplication by a rational function. So the Picard-Fuchs operators
\[\mP (\mathcal{X},\omega)\;\;\textup{and}\;\;\; \mP(\mathcal{X'},\omega')\]
are strictly equivalent.\hfill $\Diamond$.

\subsection{Moduli spaces}
We consider a smooth Calabi-Yau threefold $X$ with $h^{12}=1$. By the
famous theorem of {\sc Bogomolov}, {\sc Tian} and {\sc Todorov}, the
local deformation theory of any Calabi-Yau manifold is unobstructed,
and so by the classical deformation theory of {\sc Kodaira}, {\sc
  Spencer} \cite{KS} and {\sc Kuranishi} \cite{Kur} $X$ posses a
versal deformation over a smooth $1$-dimensional disc $\Delta$: 
\[ 
\begin{diagram}
\node{X} \arrow{e,J} \arrow{s} \node{\mathcal{X}} \arrow{s,r}{f}\\
\node{ \{0\} } \arrow{e,J} \node{\Delta}
\end{diagram}
\]
Even if $X$ is projective there will not exist a well-defined moduli space for $X$, but rather one has a moduli stack.
However, when we choose an ample line bundle $L$ on $X$, we can form families over corresponding quasi 
projective moduli spaces,\cite{Vie}. As a result, we may produce many apriori different projective families
\[f: \mathcal{X}_L \lra S_L\] 
which all have $X$ as fibre.  Now versality of the Kuranishi family $\mathcal{X} \lra \Delta$ implies that if $X$ appears 
as fibre $f^{-1}(s)$ of any projective family $f: \mathcal{X'} \lra S$ over a
curve $S$, and having a non-zero Kodaira-Spencer map at $s$, then this family is locally analytic isomorphic to 
the above model family $\mathcal{X} \lra \Delta$.\\ 

Let us call two projective families $f: \mathcal{X} \lra S$ and 
$f': \mathcal{X}' \lra S'$ {\em related},
if there exists curve $D$ and finite maps $g:D \lra S$, $g':D \lra S'$ 
and an isomorphism 
\[ \varphi: g^*\mathcal{X} \lra g'^*\mathcal{X}'\]

Clearly, if the families $\mathcal{X} \lra S$ and $\mathcal{X} \lra S'$
are both obtained as pull-back from a single family over $T$, then the 
families are related.\\

{\bf Theorem 3:} If  $\mathcal{X} \lra S$ and  $\mathcal{X}' \lra S'$
are families of Calabi-Yau manifolds, $\omega$ and $\omega'$ volume forms on 
$\mathcal{X}$ and $\mathcal{X}'$, then the Picard-Fuchs operators 
$\mP(\mathcal{X},\omega)$ and $\mP(\mathcal{X}',\omega')$
are related.\\

{\bf Theorem 2:} Let $f_1:\mathcal{X}_1 \lra S_1$ and $f_2: \mathcal{X}_2 \lra S_2$ be two projective families and
\[ \hat{\varphi}: \hat{X}_1 \lra \hat{X}_2 \]
an isomorphism between between the formal neighbourhoods $\hat{X}_i$ of a fibre $X_i=f_i^{-1}(s_i) \subset \mathcal{X}_i$,
then there exists an isomorphisms $\phi$, $\psi$
\[
\begin{diagram}
\node{U_1} \arrow{e,t}{\phi} \arrow{s,l}{f_1} \node{U_2} \arrow{s,r}{f_2}\\
\node{ V_1 } \arrow{e,t}{\psi} \node{V_2}
\end{diagram}
\]

of \'etale neighbourhoods $U_i$ of $X_i \subset \mathcal{X}_i$, $s_i
\in S_i$ ($i=1,2$) and thus the families $\mathcal{X}_1 \lra S_1$ and
$\mathcal{X}_2 \lra S_2$ are related.\\ 

{\bf proof:} This is a particular case of a very general Theorem (1.7) (Uniqueness) proven by {\sc Artin} in \cite{Art}. It states that if $F$ is a functor of locally of finite presentation and $\overline{\xi} \in F(\overline{A})$ an effective versal family, then the triple $(X,x,\xi)$ is unique up to local isomorphism for the \'etale topology, meaning that if $(X',x',\xi')$ is another algebraisation, then there is a third one, dominating both.\\

A more down to earth proof can be given by an application of the 
{\em nested approximation theorem},\cite{KPPRM}, \cite{Ron}, theorem 5.2.1.
In the above situation it can be applied as follows: assume we have
two algebraisations $\mathcal{X} \lra S$, $\mathcal{X'} \lra S'$,
given by equations of  
the form $F(x,t)=0, G(x',t')=0$.
We are looking for a algebraic maps 
\[ \phi(x,t),\;\;\;\psi(t)\]
that map $\mathcal{X}$ to $\mathcal{X}'$ and $S$ to $S'$. 
Hence we look for solutions to the system of equations
\[ F(x,t)=0,\;\;\;G(\phi(x,t),\psi(t))=0\]
Now as $\mathcal{X}$ and $\mathcal{X}'$ are analytically equivalent, we know that there exist formal solutions (or even analytic) $\phi, \psi$ to the above equations and hence by the nested approximation theorem we obtain a solution in the ring of algebraic power series.
\hfill $\Diamond$.\\

{\bf Theorem 3:} If  $\mathcal{X} \lra S$ and  $\mathcal{X}' \lra S'$
are related families of Calabi-Yau manifolds, $\omega$ and $\omega'$ 
volume forms on 
$\mathcal{X}$ and $\mathcal{X}'$, then the Picard-Fuchs operators 
$\mP(\mathcal{X},\omega)$ and $\mP(\mathcal{X}',\omega')$
are related.\\

{\bf Corollary:} If $X$ is Calabi-Yau threefold with $h^{12}=1$, which appears
as fibre in any two families $\mathcal{X}_i \lra S_i$ ($i=1,2$). Assume that
the Kodaira-Spencer map of both families is non-zero at $X$. Then the
Picard-Fuchs operator $\mathcal{P}_1$ for $\mathcal{X}_1$ and 
$\mathcal{P}_2$ for $\mathcal{X}_2$ are related.\\

We claim that the following theorem holds, but the details will
be given elsewhere.\\

{\bf Theorem 4:}\\
If $\mathcal{X}_1 \lra S_1$ and $\mathcal{X}_2 \lra S_2$ are two smooth
families of Calabi-Yau threefolds. Let 
$X_1 \subset \mathcal{X}_1$ and $X_2 \subset \mathcal{X}_2$ be
two fibres. Assume that\\
1) the Kodaira-Spencer map at $X_1$ and $X_2$ is non-zero.\\
2) $X_1$ and $X_2$ are birational.\\
Then the two families are related.\\

Theorem ${\bf 4}$ would be very useful. For example, the Calabi-Yau 
threefolds of arrangements ${\bf 152}$ and ${\bf 153}$ both have $h^{11}=41$, so a priori could be birational. But no smooth fibre of Arr. 152/198 is 
birational to a smooth fibre of Arr.153/197. The reason is, that
the Picard-Fuchs operator of ${\bf 152}$ has $KCCC$ as singularities, 
whereas that of ${\bf 153}$ has $ACCC$ singularities, so these operators 
are not related. We do not know of any other means of distinguishing members of these two families. We see that the Picard-Fuchs operator can be used as a powerful new birational invariant of a variety. We expect this to be applicable in much greater generality.\\

\section{Concluding remarks and questions}
As for families of $K3$-surfaces or other varieties, there is nothing special about families 
of Calabi-Yau three-folds that avoid having degenerations with maximal unipotent monodromy points 
and still have monodromy Zariski dense in $Sp_4(\C)$. One may ask further questions about the 
possible distribution of types of singularities. For example, the examples with second order Picard-Fuchs
operators give variations which have only $K$-points. Do there exist families with only $K$-points and  Picard-Fuchs operator of order four? Similarly, do there exists families with only $C$-points? It seems that only the current general lack of examples is the reason for this kind of ignorance.\\

As for mirror symmetry, the SYZ-approach via dual torus fibrations (T-duality), \cite{SYZ}, does in no
way presupposes the presence of a MUM-point in the moduli space. The point is rather that the
absence of a MUM-point does not give a clue where to look for an appropriate torus fibration.
However, there is a the well-known idea, going back to {\sc Miles Reid} \cite{Rei}, that the different families  of Calabi-Yau threefolds may all be connected via geometrical transitions. The most prominent 
such transition is the {\em conifold transition}, where we contract some lines on a Calabi-Yau threefold 
to form nodes,  and smooth these out to produce another Calabi-Yau threefold with different Hodge numbers. As has been suggested by {\sc Morrison},\cite{Mor2}, mirror symmetry can, in some sense, be prolongated over such transitions, and in the case of the above mentioned conifold transition, the mirror symmetric process is that of nodal degeneration, followed by a 
crepant resolution. So if we connect an orphan family to another family that has a MUM-point via a transition, one is tempted to try construct a mirror manifold for orphans using this transition to a family with a MUM-point.\\

We have seen in section 1 that the rigid Calabi-Yau from arrangement ${\bf 69}$ appears as a member of the 
orphan family ${\bf 70}$. But ${\bf 69}$ also is member of the family ${\bf 100}$, which contains even
two MUM-points.

\begin{center}
\includegraphics[height=5cm]{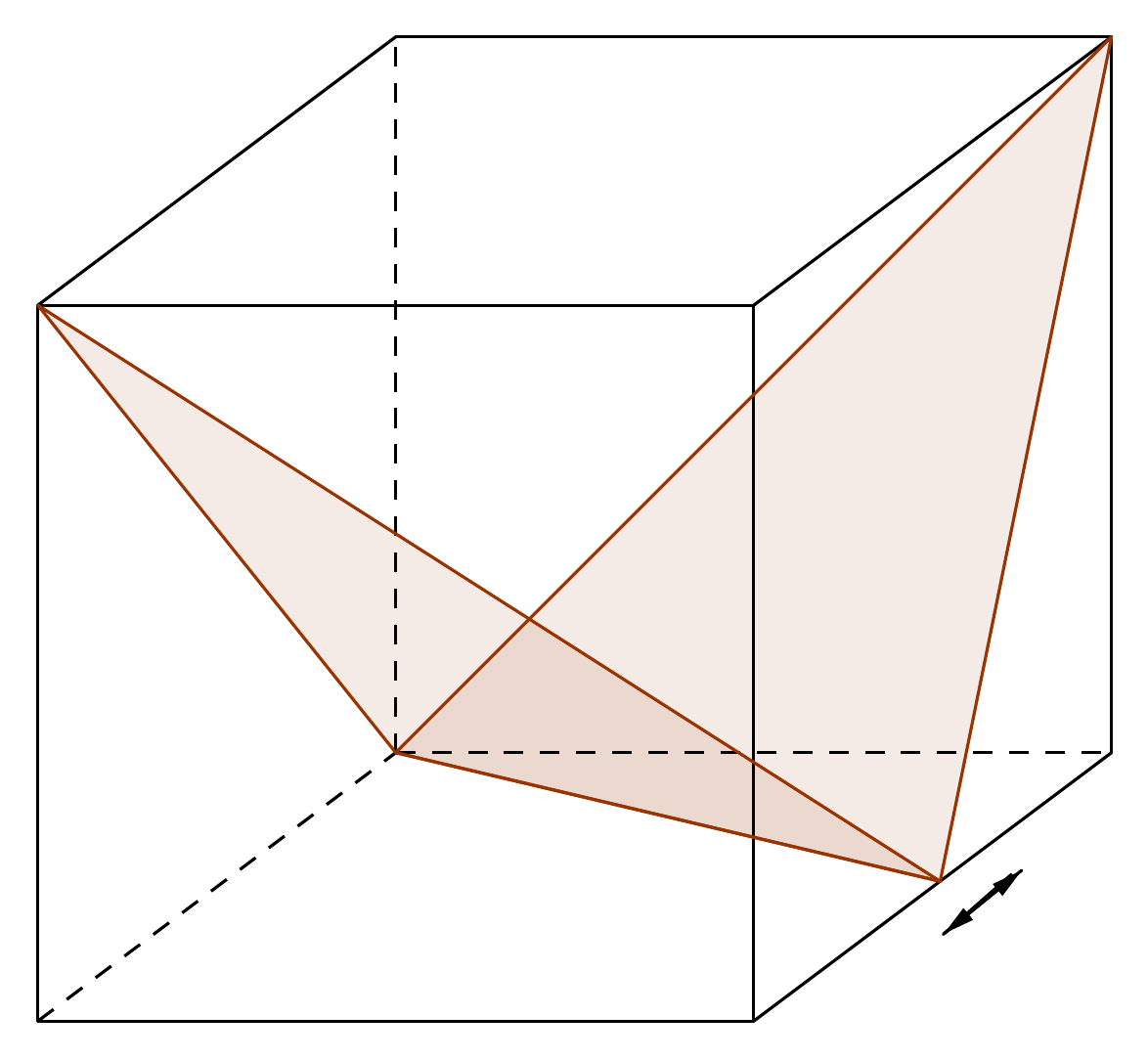}

{\bf \em Family of arrangements 100}
\end{center}

So we have the following picture

\begin{center}
\includegraphics[height=3cm]{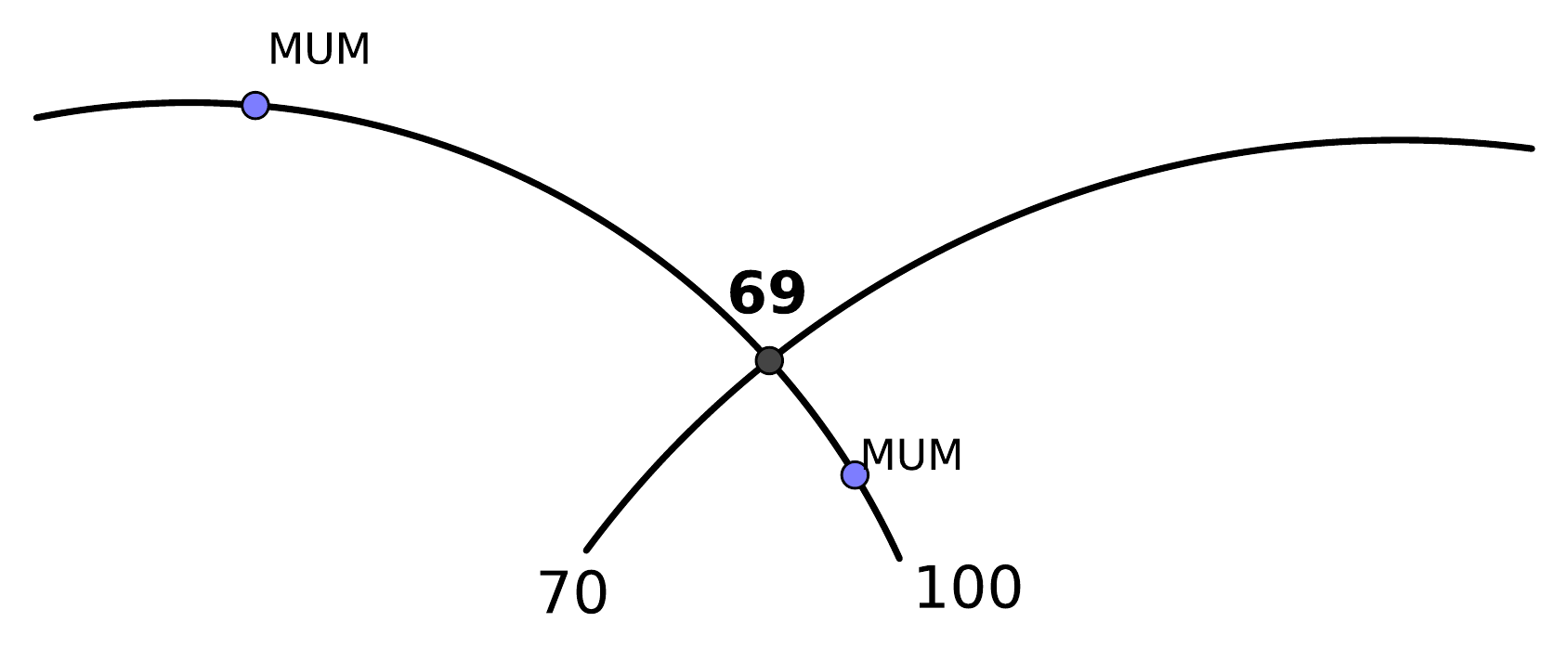}
\end{center}

In the second part of the paper we will report on the double octics families that contain a MUM-point.\\

\section{Appendix A}

We will encounter certain modular forms over and over again. For the convenience of the reader we list here the first few Fourier coefficients of these forms.\\

{\em Weight two modular form:}
Associated to the elliptic curve 
\[ y^2=x^3-x\]
of conductor $2^5=32$ there is the unique cusp form, \cite{HoMi}:
\[f_{32} :=q\prod_{n=1}^{\infty} (1-q^{4n})^2(1-q^{8n})^2 \in S_2(\Gamma_0(32))\]

\[
\begin{array}{|c||c|c|c|c|c|c|c|c|c|c|}
\hline
Name&2&3&5&7&11&13&17&19&23&29\\
\hline
\hline
f_{32}&0&0&-2&0&0&6&2&0&0&-10\\
\hline
\end{array}
\]

{\em Weight $3$ modular forms:} The weight three cusp forms for $\Gamma_0(N)$ all will appear with a non-trivial character that will play no role for us. The forms of level $8$ and $16$ are uniquely determined by there level,
which we use to name them.
\[
\begin{array}{|r||r|r|r|r|r|r|r|r|r|r|}
\hline
Name&2&3&5&7&11&13&17&19&23&CM-type\\
\hline
\hline
16&0&0&-6&0&0&10&-30&0&0&\Q(\sqrt{-1})\\
\hline
8&-2&-2&0&0&14&0&2&-34&0&\Q(\sqrt{-2})  \\
\hline
\end{array}
\]
These two forms are $\eta$-products:
\[ 8:=q\prod_{n=1}^{\infty}(1-q^n)^2(1-q^{2n})(1-q^{4n})(1-q^{8n})^2\]
\[ 16:=q\prod_{n=1}^{\infty}(1-q^{4n})^6\]
The Galois representation associated to the form 8 is the tensor square
of the Galois action associated to the form $f_{32}$. 
For more information on weight $3$ forms we refer to \cite{Schu1}.\\

{\em Weight $4$ cusp forms}\\

For the weight $4$ cusp forms for $\Gamma_0(N)$ we will used the notation used in the book of {\sc Meyer} \cite{Mey}.

\[
\begin{array}{|r||r|r|r|r|r|r|r|r|r|}
\hline
Name&2&3&5&7&11&13&17&19&23\\
\hline
\hline
6/1 &-2&-3&6&-16&12&38&-126&20&168\\
\hline
8/1 &0&-4&-2&24&-44&22&50&44 &-56\\
\hline
12/1& 0& 3&-18& 8& 36& -10& 18& -100& 72\\
\hline
32/1&0&0& 22& 0&  0&-18&-94&0 &0\\
\hline
32/2&0&8&-10&16&-40&-50&-30&40&48\\
\hline
\end{array}
\]

The first two of these forms are also nice $\eta$-products:
\[6/1= q\prod_{n=1}^{\infty}(1-q^n)^2(1-q^{2n})^2(1-q^{3n})^2(1-q^{6n})^2\]
\[8/1= q\prod_{n=1}^{\infty}(1-q^{2n})^4(1-q^{4n})^4\]
The Galois representation associated to the form 32/1 is the tensor cube
of the Galois action associated to the form $f_{32}$.

There is also one Hilbert modular form $h$ for the field $\Q(\sqrt{2})$ that plays a role.
It is the unique cuspform of weight $(4,2)$ and level $6\sqrt{2}$; its coefficients for the
first inert primes are:
\[
\begin{array}{|c||c|c|c|c|c|c|c|c|c|}
\hline
Name&2&3&5&7&11&13&17&19&23\\
\hline
h&0&9&10&4\sqrt{2}+16&-726&2938&16\sqrt{2}-62&6650&-8\sqrt{2}+40\\
\hline
\end{array}  
\]

\subsection{Appendix B}
\rule{0pt}{1pt}\medskip

\def\tfrac#1#2{{\textstyle \frac#1#2}}
\def\rsline#1 #2 #3 #4 #5 {#1&#2&#3&#4&#5} 
\def\rru{\rule{16mm}{0mm}}
\def\rrr{\rule{18mm}{0mm}}
\def\rrru{\rule{36mm}{0mm}}
\def\bb{\)\\\rrr\(}
\def\opn#1 {\textbf{#1. \quad}}
\def\pfo{\\[1mm]
\rru}
\long\def\rsym{\par\rrru }

\bgroup
\tiny

\opn4
\(xyzv  \left( x+y \right)  \left( x+ty+tz-v \right)  \left( x+y+tz-v \right) \left( y+z \right) \)
\pfo \({\Theta}^{2}-t \left( \Theta+\frac12 \right) ^{2}
\)

 \rsym
 \(\left\{
 \begin{array}{ccc}
   1&0&\infty\\\hline
   0&0&1/2\\0&0&1/2
 \end{array}\right\}
\)

\opn13
\(xyzv \left( z+y \right)  \left( x-z-v \right)  \left( x+y \right)  \left( x-z+tv \right) 
\)
\pfo\(\Theta^{2}+t(\Theta+\frac12)^{2}\)

 \rsym
 \(\left\{
 \begin{array}{ccc}
   0&-1&\infty\\\hline
   0&0&1/2\\0&0&1/2
 \end{array}\right\}
\)

\opn34
\(xyzv \left( x+y \right)  \left( x+z \right)  \left( x+y+z+v \right)  \left( y-z+tv \right) 
\)
\pfo
\(\Theta^{2}-t^{2}(\Theta+\frac12)^{2}\)

 \rsym
 \(\left\{
 \begin{array}{cccc}
   1&0&-1&\infty\\\hline
   0&0&0&1/2\\1/2&0&1/2&1/2
 \end{array}\right\}
\)

\opn72
\(xyzv  \left( x-y-v \right)  \left( x+y+z \right)  \left( y+tz+tv
\right) \left( y+z+v \right)  
\)
\pfo
\({\Theta}^{2}
+t \left( -3\,{\Theta}^{2}-2\,\Theta-\frac12 \right) 
+{t}^{2} \left( \Theta+1 \right)  \left( 2\,\Theta+1 \right) 
\)

 \rsym
 \(\left\{
 \begin{array}{cccc}
   1&1/2&0&\infty\\\hline
   0&0&0&1/2\\0&1/2&0&1
 \end{array}\right\}
\)

\opn261
\(xyzv \left( x-y-z+v \right)  \left( x+y+z+v \right)  \left( x-y+tz-tv \right)  \left( x+y+tz+tv \right) 
\)
\pfo
\({\Theta}^{2}-{t}^{2} \left( \Theta+1 \right) ^{2}\)

 \rsym
 \(\left\{
 \begin{array}{cccc}
   1&0&-1&\infty\\\hline
   0&0&0&1\\0&0&0&1
 \end{array}\right\}
\)

\opn264
\(xyzv   \left( x+(2-t)v+2\,y-(2-t)\,z \right)
\left( -x-y+2\,z-(2-t)v \right)  \left( x+y+tz \right)  
\left( y+2-2\,z \right)\)
\pfo
\({\Theta}^{2}-1/4\,{t}^{2} \left( \Theta+1 \right) ^{2}
\)

 \rsym
 \(\left\{
 \begin{array}{cccc}
   2&0&-2&\infty\\\hline
   0&0&0&1\\0&0&0&1
 \end{array}\right\}
\)

\opn270
\(xyzv  \left( x+y+z \right) \left( y+z+v \right)   \left(
  xt-2\,y+tz+tv \right)  \left( -x -2\,y+tz-v \right) 
\)
\pfo
\({\Theta}^{2}
+t \left( \frac32\,{\Theta}^{2}+\frac32\,\Theta+\frac12 \right) 
+\frac12\,{t}^{2} \left( \Theta+1 \right) ^{2}
\)

 \rsym
 \(\left\{
 \begin{array}{cccc}
   0&-1&-2&\infty\\\hline
   0&0&0&1\\0&0&0&1
 \end{array}\right\}
\)

\egroup

\subsection{Appendix C}
\rule{0pt}{1pt}\medskip

\def\tfrac#1#2{{\textstyle \frac#1#2}}
\def\rsline#1 #2 #3 #4 #5 {#1&#2&#3&#4&#5} 
\def\rru{\rule{16mm}{0mm}}
\def\rrr{\rule{18mm}{0mm}}
\def\rrru{\rule{36mm}{0mm}}
\def\bb{\)\\\rrr\(}
\def\opn#1 {\textbf{#1. \quad}}
\def\pfo{\\[1mm]
\rru}
\long\def\rsym{\par\rrru }

\bgroup
\tiny

  \textbf{33.\qquad }
\( xyzv  \left( x+y \right)\left( y+z \right)  \left( x-z+v \right)   \left( x-y-z+tv \right) 
\)
\\[1mm]
\rru\(\displaystyle {\Theta}^{2} \left( \Theta-1 \right) ^{2}
-\frac1{8}\,t{\Theta}^{2} \left( 20\,{\Theta}^{2}+3 \right) 
+\frac1{16}\,{t}^{2} \left( 8\,{\Theta}^{2}+8\,\Theta+3 \right)  \left( 2\,\Theta+1 \right) ^{2}
-\frac{1}{32}\,{t}^{3} \left( 2\,\Theta+3 \right) ^{2} \left( 2\,\Theta+1\right) ^{2}
\)

\rrru
\(\displaystyle
  \left\{ \begin {array}{cccc}   0 & 1&   2& \infty \\\hline
              0&0&0&1/2\\
              0&1/2&1&1/2\\
              1&1/2&1&3/2\\
              1&1&2&3/2
\end {array} \right\} \)

\medskip

\textbf{35.\quad }
\(xyz \left( x-v \right)\left( y-v \right) \left( z-v \right)   \left( x-y \right)
  \left( x+ty+(1-t)z-v \right)
\)
\\[1mm]
\rru
\(\Theta\, \left( \Theta-1 \right)  \left(\Theta-\frac12 \right) ^{2}
-\frac14\,t{\Theta}^{2} \left( 4\,{\Theta}^{2}+3 \right) 
-\frac14\,{t}^{2} \left( {\Theta}^{2}+\Theta+1 \right)  \left( 2\,\Theta+1 \right) ^{2}
+\frac14\,{t}^{3} \left( 2\,\Theta+3 \right)  \left( 2\,\Theta+1 \right)  \left( \Theta+1 \right) ^{2}
\)

\rrru\(\left\{ \begin {array}{cccc}1&0&-1&\infty  \\
                 \hline 0&0&0&1/2\\
                 0&1/2&1&1\\
                 1&1/2&1&1\\
                 1&1&2&3/2
            \end {array} \right\} \)

\medskip

\textbf{70.\quad}
\(
yxzv  \left( x+ty \right)\left( y-z-v \right)  \left( x-y-v \right)  \left( x-y+z \right)  
\)\\[1mm]
\rru\(
{\Theta}^{2} \left( \Theta-1 \right) ^{2}
+\frac18\,t{\Theta}^{2} \left( 20\,{\Theta}^{2}+3 \right) 
+\frac1{16}\,{t}^{2} \left( 8\,{\Theta}^{2}+8\,\Theta+3 \right)  \left( 2\,\Theta+1 \right) ^{2}
+\frac1{32}\,{t}^{3} \left( 2\,\Theta+3 \right) ^{2} \left( 2\,\Theta+1 \right) ^{2}
\)

\rrru
\(\left\{ \begin {array}{cccc} 0 & -1& -2& \infty  \\
            \hline 0&0&0&1/2\\
            0&1/2&1&1/2\\
            1&1/2&1&3/2\\
            1&1&2&3/2
\end {array} \right\} \)

\medskip

\textbf{71.\quad}
\(xyzv \left( x+y \right)  \left( x+y+z+v \right)  \left( ty-z-v \right)  \left( -x+ty-z \right) 
\)\\[1mm]
\rru\(
\Theta\, \left( \Theta-1 \right)  \left( \Theta -\frac12 \right) ^{2}
+t{\Theta}^{2} \left( 4\,{\Theta}^{2}+1 \right) 
+\frac1{16}\,{t}^{2} \left( 20\,{\Theta}^{2}+20\,\Theta+9 \right)  \left( 2\,\Theta+1 \right) ^{2}
+\frac18\,{t}^{3} \left( 2\,\Theta+3 \right) ^{2} \left( 2\,\Theta+1 \right) ^{2}
\)

\rrru
\(\left\{
  \begin {array}{cccc}
    0&-1/2&-1&\infty\\\hline
    0&0&0&1/2\\
    1/2&1&1/2&1/2\\
    1/2&1&1/2&3/2\\
    1&2&1&3/2
\end {array} \right\} \)

\medskip

\textbf{97.\quad }
\(xyzv   \left( x+y \right)  \left( x+y+z+v \right)  \left( -x+tz-v \right) \left( y-v-z \right)
\)\\[1mm]
\rru
\(\Theta\, \left( \Theta-1 \right)  \left( \Theta -\frac 12\right) ^{2}
+\frac18\,t{\Theta}^{2} \left( 20\,{\Theta}^{2}+7 \right) \
+\frac1{32}\,{t}^{2} \left( 16\,{\Theta}^{2}+16\,\Theta+7 \right)  \left( 2\,\Theta+1 \right) ^{2}
+\frac1{32}\,{t}^{3} \left( 2\,\Theta+3 \right) ^{2} \left( 2\,\Theta+1 \right) ^{2}
\)

\rrru
\(\left\{
  \begin {array}{cccc}0&-1&-2&\infty  \\
    \hline
    0&0&0&1/2\\
    1/2&0&1&1/2\\
    1/2&1&1&3/2\\
    1&1&2&3/2
    \end {array} \right\} \)

\medskip
  
\textbf{98.\quad }
\(xyzv   \left( x+z-v \right)  \left( x+z+y \right) \left( y+z+v \right) \left( y+tz+tv \right) 
\)\\[1mm]
\rru
\(
{\Theta}^{2} \left( \Theta-1 \right) ^{2}
-\frac14\,t{\Theta}^{2} \left( 16\,{\Theta}^{2}+3 \right) 
+\frac14\,{t}^{2} \left( 5\,{\Theta}^{2}+5\,\Theta+3 \right)  \left( 2\,\Theta+1 \right) ^{2}
-\frac12\,{t}^{3} \left( 2\,\Theta+3 \right)  \left( 2\,\Theta+1 \right)  \left( \Theta+1 \right) ^{2}
\)

\rrru
\(\left\{
  \begin {array}{cccc}1&1/2&0&\infty \\
    \hline
    0&0&0&1/2\\
    0&1&0&1\\1&1&1&1\\
    1&2&1&3/2
\end {array} \right\} \)

\medskip

\opn152
\(
xyzv   \left( x+v-y-z \right)  \left( x+y+z+v \right)  \left(
  x-y+tz-tv \right) \left( y+v \right) 
\)
\pfo
\(
\Theta\, \left( \Theta-1 \right)  \left( \Theta -\frac 12 \right) ^{2}
+\frac12\,t\Theta\, \left( 2\,{\Theta}^{3}-8\,{\Theta}^{2}+6\,\Theta-1 \right) 
+{t}^{2} \left( -2\,{\Theta}^{4}-4\,{\Theta}^{3}-\frac{11}4\,{\Theta}^{2}-{\frac {17}{4}}\,\Theta-{\frac {11}{16}} \right) \bb
+{t}^{3} \left( -2\,{\Theta}^{4}+\frac14\,{\Theta}^{2}+\frac72\,\Theta+{\frac {9}{8}} \right) 
+\frac1{16}\,{t}^{4} \left( 2\,\Theta+1 \right)  \left( 8\,{\Theta}^{3}+44\,{\Theta}^{2}+62\,\Theta+25 \right) 
\bb+\frac14\,{t}^{5} \left( 2\,\Theta+3 \right)  \left( 2\,\Theta+1 \right)  \left( \Theta+1 \right) ^{2}
\)
\rsym
\(\left\{
  \begin {array}{cccc}
    1&0&-1&\infty\\\hline
    0&0&0&1/2\\
    1/2&1/2&0&1\\
    1/2&1/2&2&1\\
    1&1&2&3/2
\end {array} \right\} \)

\medskip
\opn153
\(xyzv   \left( x+z+y \right)  \left( -x+ty-v \right)  \left( -x+ty-z-v \right) \left( y+z+v \right)
\)
\pfo
\(\Theta\, \left( \Theta-1 \right)  \left(\Theta -\frac12\right) ^{2}
+\frac18\,t\Theta\, \left( 28\,{\Theta}^{3}-16\,{\Theta}^{2}+17\,\Theta-2 \right) 
+{t}^{2} \left( {\frac {19}{4}}\,{\Theta}^{4}+\frac72\,{\Theta}^{3}+{\frac {39}{8}}\,{\Theta}^{2}+{\frac {13}{8}}\,\Theta+{\frac {19}{64}} \right) \bb
+{t}^{3} \left( {\frac {25}{8}}\,{\Theta}^{4}+6\,{\Theta}^{3}+{\frac {109}{16}}\,{\Theta}^{2}+\frac72\,\Theta+{\frac {89}{128}} \right) 
+{\frac {1}{64}}\,{t}^{4} \left( 2\,\Theta+1 \right)  \left( 32\,{\Theta}^{3}+80\,{\Theta}^{2}+82\,\Theta+29 \right) 
\bb+\frac1{32}\,{t}^{5} \left( 2\,\Theta+3 \right)  \left( 2\,\Theta+1 \right)  \left( \Theta+1 \right) ^{2}
\)
\rsym
\(\left\{
  \begin {array}{cccc}0&-1&-2&\infty  \\\hline
    0&0&0&1/2\\
    1/2&1/2&1/2&1\\
    1/2&1/2&3/2&1\\
    1&1&2&3/2
\end {array} \right\} \)

\medskip
\opn197
\(xyzv \left( x-y-z+v \right) \left( x+tz+v \right)  \left( x+ty+tz \right) \left( ty+tz+v \right)   
\)
\pfo
\({\Theta}^{2} \left(\Theta - \frac12\right)  \left( \Theta+\frac12 \right) 
+\frac18\,t \left( 2\,\Theta+1 \right)  \left( 32\,{\Theta}^{3}+16\,{\Theta}^{2}+18\,\Theta+5 \right) 
+{t}^{2} \left( 25\,{\Theta}^{4}+52\,{\Theta}^{3}+{\frac {121}{2}}\,{\Theta}^{2}+37\,\Theta+{\frac {145}{16}} \right) \bb
+{t}^{3} \left( 38\,{\Theta}^{4}+124\,{\Theta}^{3}+183\,{\Theta}^{2}+133\,\Theta+{\frac {307}{8}} \right) 
+{t}^{4} \left( \Theta+1 \right)  \left( 28\,{\Theta}^{3}+100\,{\Theta}^{2}+133\,\Theta+63 \right) 
\bb+2\,{t}^{5} \left( \Theta+2 \right)  \left( \Theta+1 \right)  \left( 2\,\Theta+3 \right) ^{2}
\)
\rsym
\(\left\{
  \begin {array}{cccc}0&-1/2&-1&\infty  \\\hline
    -1/2&0&0&1\\
    0&1/2&1/2&3/2\\
    0&3/2&1/2&3/2\\
    1/2&2&1&2
\end {array} \right\} \)

\medskip
\opn198
\(xyzv  \left( x-y-v \right)  \left( x+y+z \right)  \left( x-y-z+tx \right) \left( y+z+v \right) 
\)
\pfo
\(\Theta\, \left( \Theta-1 \right)  \left( \Theta-\frac12 \right) ^{2}
+\frac18\,t{\Theta}^{2} \left( 24\,{\Theta}^{2}+5 \right) 
+{t}^{2} \left( {\frac {13}{4}}\,{\Theta}^{4}+\frac{13}2\,{\Theta}^{3}+{\frac {81}{16}}\,{\Theta}^{2}+{\frac {29}{16}}\,\Theta+{\frac {5}{16}} \right) \bb
+{t}^{3} \left( \frac32\,{\Theta}^{4}+6\,{\Theta}^{3}+8\,{\Theta}^{2}+4\,\Theta+{\frac {25}{32}} \right) 
+{\frac {1}{64}}\,{t}^{4} \left( 2\,\Theta+5 \right) ^{2} \left( 2\,\Theta+1 \right) ^{2}
\)
\rsym
\(\left\{
  \begin {array}{cccc}0&-1&-2&\infty  \\\hline
    0&0&0&1/2\\
    1/2&1/2&1/2&1/2\\
    1/2&1/2&1/2&5/2\\
    1&1&1&5/2
\end {array} \right\} \)

\medskip
\opn243
\(xyzv   \left( x+y+v \right)  \left( x+y+z \right)  \left( x+ty+z+v \right) \left( y+z+v \right)
\)
\pfo
\(\Theta\, \left( \Theta-2 \right)  \left( \Theta-1 \right) ^{2}
-\frac16\,t\Theta\, \left( \Theta-1 \right)  \left( 19\,{\Theta}^{2}-19\,\Theta+9 \right) 
+\frac13\,{t}^{2}{\Theta}^{2} \left( 11\,{\Theta}^{2}+4 \right) 
-\frac1{24}\,{t}^{3} \left( 11\,{\Theta}^{2}+11\,\Theta+5 \right)
\left( 2\,\Theta+1 \right) ^{2} \bb
+\frac1{48}\,{t}^{4} \left( 2\,\Theta+3 \right) ^{2} \left( 2\,\Theta+1 \right) ^{2}
\)
\rsym
\(\left\{
  \begin {array}{ccccc}2&3/2&1&0&\infty\\\hline
    0&0&0&0&1/2\\
    1&1&1/2&1&1/2\\
    1&1&1/2&1&3/2\\
    2&2&1&2&3/2
\end {array} \right\} \)

\medskip
\opn247
\(xyzv   \left( x-y-v \right)  \left( x+y+z \right)  \left( -x+tz-tv \right) \left( y+z+v \right)
\)
\pfo
\(
{\Theta}^{2} \left( \Theta-1 \right) ^{2}
+t{\Theta}^{2} \left( 5\,{\Theta}^{2}+1 \right) 
+{t}^{2} \left( 2\,{\Theta}^{2}+2\,\Theta+1 \right)  \left(
  2\,\Theta+1 \right) ^{2} 
+{t}^{3} \left( 2\,\Theta+3 \right)  \left( 2\,\Theta+1 \right)  \left( \Theta+1 \right) ^{2}
\)
\rsym
\(\left\{
  \begin {array}{cccc}0&-1/2&-1&\infty  \\\hline
    0&0&0&1/2\\
    0&1/2&1&1\\
    1&1/2&1&1\\
    1&1&2&3/2
\end {array} \right\} \)

\medskip
\opn248
\(xyzv   \left( x+z+v \right)  \left( x+y+z \right)  \left( x+(y+1)y-tz+v \right) \left( y-z-v \right)
\)
\pfo
\(\Theta\, \left( \Theta-2 \right)  \left( \Theta-1 \right) ^{2}
+\frac16\,t\Theta\, \left( \Theta-1 \right)  \left( 37\,{\Theta}^{2}-61\,\Theta+36 \right)
+\frac16\,{t}^{2}\Theta\, \left( 91\,{\Theta}^{3}-124\,{\Theta}^{2}+121\,\Theta-36 \right)  \bb
+{t}^{3} \left( {\frac {115}{6}}\,{\Theta}^{4}-\frac53\,{\Theta}^{3}+{\frac {107}{6}}\,{\Theta}^{2}+\frac23\,\Theta+\frac12 \right) 
+{t}^{4} \left( {\frac {79}{6}}\,{\Theta}^{4}+16\,{\Theta}^{3}+{\frac
    {113}{6}}\,{\Theta}^{2}+8\,\Theta+\frac32 \right)\bb 
+\frac16\,{t}^{5} \left( 2\,\Theta+1 \right)  \left(
  14\,{\Theta}^{3}+29\,{\Theta}^{2}+27\,\Theta+9 \right) 
+\frac16\,{t}^{6} \left( 2\,\Theta+3 \right)  \left( 2\,\Theta+1
\right)  \left( \Theta+1 \right) ^{2} 
\)
\rsym
\(\left\{
  \begin {array}{cccccc} 0&-1/2&-1&-3/2&-2&\infty\\\hline
    0&0&0&0&0&1/2\\
    1&1&0&1&1&1\\
    1&1&2&1&1&1\\
    2&2&2&2&2&3/2
\end {array} \right\} \)

\medskip
\opn250
\(xyzv   \left( x+y+z \right)  \left( x+ty-z+v \right) \left( x+z+v \right)  \left( y+z-v \right)
\)
\pfo
\(\Theta\, \left( \Theta-1 \right)  \left( \Theta -\frac12\right) ^{2}
+\frac18\,t\Theta\, \left( 44\,{\Theta}^{3}-96\,{\Theta}^{2}+65\,\Theta-12 \right) 
+{t}^{2} \left( \frac{19}2\,{\Theta}^{4}-23\,{\Theta}^{3}+{\frac {131}{8}}\,{\Theta}^{2}-{\frac {47}{8}}\,\Theta-\frac14 \right) 
\bb+{t}^{3} \left( \frac52\,{\Theta}^{4}-20\,{\Theta}^{3}-{\frac
    {23}{4}}\,\Theta-{\frac {17}{32}} \right) 
-\frac1{32}\,{t}^{4} \left( 68\,{\Theta}^{2}+100\,\Theta+53 \right)  \left( 2\,\Theta+1 \right) ^{2}
\bb-\frac14\,{t}^{5} \left( 8\,{\Theta}^{2}+14\,\Theta+9 \right)  \left( 2\,\Theta+1 \right) ^{2}
-\frac18\,{t}^{6} \left( 2\,\Theta+3 \right) ^{2} \left( 2\,\Theta+1 \right) ^{2}
\)
\rsym
\(\left\{
  \begin {array}{cccccc}1&0&-1/2&-1&-2&\infty  \\\hline
    0&0&0&0&0&1/2\\
    1&1/2&1&1/2&1&1/2\\
    1&1/2&3&1/2&1&3/2\\
    2&1&4&1&2&3/2
\end {array} \right\} \)

\medskip
\opn252
\(xyzv \left( x+y+v \right)  \left( x+y+z \right)  \left(-x+tz+v \right)  \left( -x-2\,y+tz-v \right) 
\)
\pfo
\({\Theta}^{2} \left( \Theta-1 \right) ^{2}
+\frac12\,t{\Theta}^{2} \left( 5\,{\Theta}^{2}+1 \right) 
+\frac14\,{t}^{2} \left( 2\,{\Theta}^{2}+2\,\Theta+1 \right)  \left( 2\,\Theta+1 \right) ^{2}
+\frac18\,{t}^{3} \left( 2\,\Theta+3 \right)  \left( 2\,\Theta+1 \right)  \left( \Theta+1 \right) ^{2}
\)
\rsym
\(\left\{
  \begin {array}{cccc} 0&-1&-2&\infty\\\hline
    0&0&0&1/2\\
    0&1/2&1&1\\
    1&1/2&1&1\\
    1&1&2&3/2
\end {array} \right\} \)

\medskip
\opn258
\(xyzv  \left( x-y+2\,z-2v \right)  \left( x-y+z-v \right)  \left( x+ty+z+tv \right) \left( y-z+2v \right) 
\)
\pfo
\(\Theta\, \left( \Theta-1 \right)  \left(\Theta -\frac12\right) ^{2}
+\frac18\,t\Theta\, \left( 20\,{\Theta}^{3}+48\,{\Theta}^{2}-21\,\Theta+6 \right) 
+{t}^{2} \left( -{\Theta}^{4}+19\,{\Theta}^{3}+{\frac {39}{2}}\,{\Theta}^{2}+{\frac {47}{8}}\,\Theta+{\frac {11}{8}} \right) \bb
+{t}^{3} \left( -5\,{\Theta}^{4}-5\,{\Theta}^{3}+{\frac {61}{2}}\,{\Theta}^{2}+{\frac {127}{8}}\,\Theta+{\frac {37}{8}} \right) 
+{t}^{4} \left( -{\Theta}^{4}-21\,{\Theta}^{3}-\frac{21}2\,{\Theta}^{2}-{\frac {9}{8}}\,\Theta+{\frac {7}{8}} \right) \bb
+\frac18\,{t}^{5} \left( 2\,\Theta+1 \right)  \left( 10\,{\Theta}^{3}-9\,{\Theta}^{2}-27\,\Theta-13 \right) 
+\frac14\,{t}^{6} \left( 2\,\Theta+3 \right)  \left( 2\,\Theta+1 \right)  \left( \Theta+1 \right) ^{2}
\)
\rsym
\(\left\{
  \begin {array}{cccccc}1&0&-1/2&-1&-2&\infty \\\hline
    0&0&0&0&0&1/2\\
    1&1/2&1&0&1&1\\
    3&1/2&1&1&1&1\\
    4&1&2&1&2&3/2
 \end {array} \right\} \)

\medskip
\opn266
\(xyzv   \left( 2x+y+2v \right)  \left( x+(t+1)y-z+v \right)  \left( x+ty+z \right) \left( y-2\,z+2v \right)
\)
\pfo
\(\Theta\, \left( \Theta-1 \right)  \left(\Theta -\frac12\right) ^{2}
+\frac14\,t\Theta\, \left( 44\,{\Theta}^{3}-48\,{\Theta}^{2}+37\,\Theta-6 \right) 
+{t}^{2} \left( 50\,{\Theta}^{4}-56\,{\Theta}^{3}+40\,{\Theta}^{2}-\frac52\,\Theta+\frac38 \right) \bb
+{t}^{3} \left( 120\,{\Theta}^{4}-288\,{\Theta}^{3}-75\,{\Theta}^{2}-105\,\Theta-21 \right) 
+{t}^{4} \left( 112\,{\Theta}^{4}-1008\,{\Theta}^{3}-718\,{\Theta}^{2}-720\,\Theta-{\frac {303}{2}} \right) \bb
+{t}^{5} \left( -224\,{\Theta}^{4}-2464\,{\Theta}^{3}-1924\,{\Theta}^{2}-1628\,\Theta-324 \right) 
+{t}^{6} \left( -960\,{\Theta}^{4}-4224\,{\Theta}^{3}-4296\,{\Theta}^{2}-2448\,\Theta-450 \right) \bb
+{t}^{7} \left( -1600\,{\Theta}^{4}-4992\,{\Theta}^{3}-6368\,{\Theta}^{2}-3504\,\Theta-696 \right)
-32\,{t}^{8} \left( 2\,\Theta+1 \right)  \left( 22\,{\Theta}^{3}+57\,{\Theta}^{2}+59\,\Theta+21 \right) \bb
-128\,{t}^{9} \left( 2\,\Theta+3 \right)  \left( 2\,\Theta+1 \right)  \left( \Theta+1 \right) ^{2}
\)
\rsym
\(\left\{
  \begin {array}{cccccccc}
    1/2&0&-1/4&-1/2&-1&(-1+\sqrt{-3})/4&(-1-\sqrt{-3})/4&\infty\\\hline
    0&0&0&0&0&0&0&1/2\\
    1&1/2&1&1/2&1&1&1&1\\
    1&1/2&1&1/2&1&3&3&1\\
    2&1&2&1&2&4&4&3/2
\end {array} \right\} \)

\medskip
\opn273
\(xyzv    \left( x+y+z \right) \left( 2\,x-2\,z-v \right)  \left( x+2\,ty-z+tv \right) \left( 2\,y+2\,z+v \right)
\)
\pfo
\(\Theta\, \left( \Theta-1 \right)  \left(\Theta - \frac12\right) ^{2}
+\frac14\,t\Theta\, \left( 44\,{\Theta}^{3}-48\,{\Theta}^{2}+37\,\Theta-6 \right) 
+{t}^{2} \left( 50\,{\Theta}^{4}-56\,{\Theta}^{3}+40\,{\Theta}^{2}-\frac52\,\Theta+\frac38 \right)\bb 
+{t}^{3} \left( 120\,{\Theta}^{4}-288\,{\Theta}^{3}-75\,{\Theta}^{2}-105\,\Theta-21 \right) 
+{t}^{4} \left( 112\,{\Theta}^{4}-1008\,{\Theta}^{3}-718\,{\Theta}^{2}-720\,\Theta-{\frac {303}{2}} \right)\bb 
+{t}^{5} \left( -224\,{\Theta}^{4}-2464\,{\Theta}^{3}-1924\,{\Theta}^{2}-1628\,\Theta-324 \right) 
+{t}^{6} \left( -960\,{\Theta}^{4}-4224\,{\Theta}^{3}-4296\,{\Theta}^{2}-2448\,\Theta-450 \right) \bb
+{t}^{7} \left( -1600\,{\Theta}^{4}-4992\,{\Theta}^{3}-6368\,{\Theta}^{2}-3504\,\Theta-696 \right) 
-32\,{t}^{8} \left( 2\,\Theta+1 \right)  \left( 22\,{\Theta}^{3}+57\,{\Theta}^{2}+59\,\Theta+21 \right) \bb
-128\,{t}^{9} \left( 2\,\Theta+3 \right)  \left( 2\,\Theta+1 \right)  \left( \Theta+1 \right) ^{2}
\)
\rsym
\(\left\{
  \begin {array}{cccccccc}
    1/2&0&-1/4&-1/2&-1&(-1+\sqrt{-3})/4&(-1-\sqrt{-3})/4&\infty\\\hline
    0&0&0&0&0&0&0&1/2\\
    1&1/2&1&1/2&1&1&1&1\\
    1&1/2&1&1/2&1&3&3&1\\
    2&1&2&1&2&4&4&3/2
\end {array} \right\} \)

\egroup

\end{document}